\documentclass[11pt,reqno]{amsart}

\usepackage[a4paper,margin=1in]{geometry}
\usepackage{amsmath,amssymb,amsthm,mathtools}
\usepackage{array,booktabs,longtable,tabularx}
\usepackage{enumitem}
\usepackage{microtype}
\usepackage[unicode,colorlinks=true,linkcolor=blue,citecolor=blue,urlcolor=blue]{hyperref}
\usepackage[nameinlink,noabbrev]{cleveref}
\newtheorem{theorem}{Theorem}[section]
\newtheorem{proposition}[theorem]{Proposition}
\newtheorem{lemma}[theorem]{Lemma}
\newtheorem{corollary}[theorem]{Corollary}
\newtheorem{criterion}[theorem]{Criterion}
\theoremstyle{definition}
\newtheorem{definition}[theorem]{Definition}
\newtheorem{example}[theorem]{Example}
\newtheorem{construction}[theorem]{Construction}
\theoremstyle{remark}
\newtheorem{remark}[theorem]{Remark}
\newtheorem{question}[theorem]{Question}

\newcommand{\ZZ}{\mathbb Z}
\newcommand{\QQ}{\mathbb Q}
\newcommand{\FF}{\mathbb F}
\newcommand{\End}{\operatorname{End}}
\newcommand{\Hom}{\operatorname{Hom}}
\newcommand{\Tor}{\operatorname{Tor}}
\newcommand{\Ext}{\operatorname{Ext}}
\newcommand{\HH}{\operatorname{HH}}
\newcommand{\im}{\operatorname{im}}
\newcommand{\Ker}{\operatorname{Ker}}
\newcommand{\cO}{\mathcal O}
\newcommand{\cA}{\mathcal A}
\newcommand{\cR}{\mathcal R}
\newcommand{\cC}{\mathcal C}
\newcommand{\Def}{\mathfrak D}
\newcommand{\UDef}{\mathfrak K}
\newcommand{\gr}{\operatorname{gr}}
\newcommand{\ann}{\operatorname{ann}}
\newcommand{\id}{\operatorname{id}}
\newcommand{\rank}{\operatorname{rank}}
\newcommand{\coker}{\operatorname{coker}}

\newcommand{\tors}{\operatorname{tors}}
\newcommand{\cS}{\mathcal S}
\newcommand{\cV}{\mathcal V}

\newcommand{\cX}{\mathcal X}
\newcommand{\cJ}{\mathcal J}

\title[Independent arithmetic defects]{Universal Operator Envelopes and Independent Arithmetic Defects of Additive Ternary $\Gamma$-Rings}

\author{Chandrasekhar Gokavarapu}
\address{Research Scholar, Department of Mathematics, Acharya Nagarjuna University, Guntur, Andhra Pradesh 522510, India}
\address{Lecturer in Mathematics, Department of Mathematics, Government College (Autonomous), Rajahmundry, Andhra Pradesh 533105, India}
\email{chandrasekhargokavarapu@gmail.com}
\email{chandrasekhargokavarapu@gcrjy.ac.in}
\thanks{Corresponding author: Chandrasekhar Gokavarapu. ORCID: \href{https://orcid.org/0009-0006-5306-371X}{0009-0006-5306-371X}.}

\author{Sajani Lavanya Madasi}
\address{Lecturer in Mathematics, Department of Mathematics, Government College (Autonomous), Rajahmundry, Andhra Pradesh 533105, India}
\email{sajanimaths@gcrjy.ac.in}

\author{Madhusudhana Rao Dasari}
\address{Research Supervisor, Department of Mathematics, Acharya Nagarjuna University, Pedakakani, Guntur, Andhra Pradesh 522510, India}
\address{Lecturer in Mathematics, Department of Mathematics, Government College for Women (Autonomous), Pattabhipuram, Guntur, Andhra Pradesh 522006, India}
\email{dmrmaths@gmail.com}


\hypersetup{
  pdftitle={Universal Operator Envelopes and Independent Arithmetic Defects of Additive Ternary Gamma-Rings},
  pdfauthor={Chandrasekhar Gokavarapu, Sajani Lavanya Madasi, and Madhusudhana Rao Dasari},
  pdfkeywords={ternary Gamma-ring, universal operator envelope, arithmetic defect, realization theorem, Smith normal form, base-change obstruction, periodic relative homology}
}

\begin{document}

\begin{abstract}
For an additive ternary $\Gamma$-ring $T$, we construct the associative rng generated by the two-element, two-index operators on ternary modules and compare it with its realization on the regular module.  For a finite-free two-step system $T_\mu=E\oplus W$, two integral matrices control complementary parts of this comparison: the middle-relation matrix $B_\mu$ presents the hidden quotient occurring in the integral regular kernel, while the visible flattening $\Theta_\mu$ controls the visible image, its saturation and its nonflat base-change correction.  The torsion cokernels of these matrices define an ordered arithmetic defect pair; this pair records selected arithmetic data and is not asserted to be the whole regular kernel.

The principal result is a constructive realization theorem: for every ordered pair $(G,H)$ of finite abelian groups, there is a torsion-free two-step ternary $\mathbb Z$-ring whose visible-saturation torsion is $G$ and whose middle-relation torsion is $H$.  Thus the two torsion mechanisms, including their prime supports and invariant factors, can be prescribed independently.  The proof builds two rank-$(2,1)$ atomic tensor families with complementary Smith profiles and establishes an orthogonal-sum theorem for their presentation matrices.  We also identify the standard base-change obstruction attached to these two matrices, describe their determinantal loci, realize classical rank-one matrix-factorization periodicity in an example not induced by an associative $\mathbb Z$-bilinear binary product on its underlying group, and classify all $256$ integral $2\times2\times2$ sign tensors by the resulting defect data.  All finite calculations use exact integer arithmetic and are independently reproducible.
\end{abstract}

\maketitle

\section{Introduction}

Ternary rings, associative triple systems and related multiary structures admit embeddings and enveloping constructions that translate ternary multiplication into ordinary associative algebra; see Lister's structural theory of ternary rings \cite{Lister1971}, Meyberg's lectures on triple systems \cite{Meyberg1972}, and the universal-envelope literature represented by \cite{Bremner2014,Elgendy2014}. The $\Gamma$-ring formalism initiated by Nobusawa and Barnes provides a different indexed multiplication mechanism \cite{Nobusawa1964,Barnes1966}. In an additive ternary $\Gamma$-module, however, the elementary operator acting on a module element depends on two elements of $T$ and two $\Gamma$-indices. A module-theoretic envelope should therefore retain relations among these four slots rather than only embed the elements of $T$ themselves.

We first give an explicit integral presentation of the operator envelope.  The existence of associative enveloping algebras representing module categories is known abstractly for algebras over operads \cite{LodayVallette2012,Khoroshkin2022}, so this presentation is foundational rather than the main novelty.  Our central problem is arithmetic: how independently can the universal operators fail to be saturated in the regular representation and fail to satisfy the universal middle relations?  For finite-free two-step systems these failures are measured by two finitely generated abelian groups
\[
\mathsf S_{\mathrm{vis}}(\mu)=\tors(\coker\Theta_\mu),
\qquad
\mathsf S_{\mathrm{mid}}(\mu)=\tors(\coker B_\mu).
\]
The main theorem proves that the ordered pair can be prescribed arbitrarily.  In particular, it is not forced by the abstract envelope, by the additive group of the ternary ring, or by either component of the pair.

The comparison between universal and regularly visible operators also has several homological layers:
\begin{enumerate}[label=\textup{(\roman*)},leftmargin=2.2em]
\item the intrinsic regular kernel, consisting of slot operators invisible on $T$;
\item the unital regular kernel created after adjoining an identity;
\item the conormal quotient $\UDef/\UDef^2$, which controls first relative homology;
\item the full graded tower $\UDef^n/\UDef^{n+1}$;
\item when $\UDef^2=0$, a Hochschild class measuring the failure of the regular quotient to split as a ring.
\end{enumerate}
We organize these standard quotient-ring invariants for the particular regular representation $\cA(T,\Gamma)\twoheadrightarrow\cR_T$ and compute them in ternary families below.

The principal results, in order of mathematical dependence, are summarized below.

\begin{enumerate}[label=\textup{(\arabic*)},leftmargin=2.2em]
\item We give an explicit structure formula for the operator envelope and regular defect of every two-step tensor system $T_\mu=E\oplus W$.  The visible flattening $\Theta_\mu$ and the middle-relation matrix $B_\mu$ determine its free defect, torsion, square-zero extension class and characteristic jumps.
\item For each $n\ge2$ we construct two rank-$(2,1)$ tensors $V_n$ and $M_n$ with complementary arithmetic profiles
\[
\mathsf S_{\mathrm{arith}}(V_n)=(\ZZ/n\ZZ,0),
\qquad
\mathsf S_{\mathrm{arith}}(M_n)=(0,\ZZ/n\ZZ).
\]
An orthogonal-sum theorem then realizes every ordered pair of finite abelian groups on a torsion-free two-step ternary ring.  This settles the realization and prime-independence questions without appealing to a binary-derived model.
\item The universal operator rng $\cO(T,\Gamma)$ has a concrete four-slot presentation, its Dorroh extension represents additive ternary $\Gamma$-modules, and this explicit envelope commutes with arbitrary scalar extension.  These statements supply the framework in which the arithmetic pair is intrinsic.
\item For the regular image $\cR_T$ and kernel $\UDef$, we apply the standard conormal, relative $\Tor$/$\Ext$, filtered-ideal and square-zero-extension constructions to this ternary regular quotient.  The ternary content lies in the canonical source of $\UDef$ and in the computations for the families below.
\item The two-matrix decomposition isolates the standard base-change correction: all new kernel failure is $\Tor^{\ZZ}_1(\coker\Theta_\mu,S)$, while $\coker B_\mu$ commutes with base change.  Maximal minors and Fitting ideals then give the corresponding bad-prime and determinantal loci.
\item A directly defined ternary ring that is not induced by an associative binary product has
\[
\cA\cong\ZZ[p]/(p^3),\qquad \cR_T\cong\ZZ[p]/(p^2).
\]
A rank-one matrix factorization gives a two-periodic free resolution and nonzero relative $\Tor_n$ and $\Ext^n$ for every $n\ge1$.
\item A nonunital rank-one family has zero intrinsic defect but a nontrivial infinite conormal tower. A second rank-one matrix factorization makes its relative homology two-periodic with
\[
\Tor_{2j+1}\cong\ZZ/d^2\ZZ,
\qquad
\Ext^{2j+2}\cong\ZZ/d^2\ZZ.
\]
\item A torsion-free rank-three ternary group defined by a sign tensor has conormal defect $\ZZ^2\oplus\ZZ/2\ZZ$. Reduction modulo $2$ creates an additional rank drop.  Applying the classical switching criterion to the middle-relation incidence matrix yields a new exact census of all $256$ sign tensors; $252$ have middle-relation torsion.
\end{enumerate}

Hochschild's interpretation of second cohomology in terms of algebra extensions \cite{Hochschild1945} is used only after the ternary-specific operator quotient has been constructed. Standard homological facts follow the conventions of \cite{CartanEilenberg1956,Weibel1994,Witherspoon2019}. All matrix computations in the finite-free tensor section use exact integer arithmetic. The accompanying notebook is intended for verification and exploration; none of the proofs relies on numerical approximation or statistical evidence.

\section{Position within the literature and precise novelty}\label{sec:literature}

The construction developed here lies at the intersection of four established theories, but it is not identical to any of them.  Since the novelty of an operator envelope can be assessed only against its nearest predecessors, we record the distinctions explicitly.

\subsection{Gamma rings and their operator rings}
Nobusawa introduced $\Gamma$-rings as indexed binary products, and Barnes isolated the formulation that became standard in ring theory \cite{Nobusawa1964,Barnes1966}.  Booth subsequently constructed left and right operator rings for a $\Gamma$-ring \cite{Booth1986}; these rings were used to transfer prime ideals, radicals, Noetherian conditions, bi-ideals and quasi-ideals between a $\Gamma$-ring and ordinary associative rings \cite{CoppageLuh1971,BoothGroenewald1992}.  Dutta and Sardar developed the corresponding operator semirings for $\Gamma$-semirings and proved ideal-lattice correspondences \cite{DuttaSardar2003}.

The slot count alone does not separate the present construction from classical
$\Gamma$-ring machinery.  Indeed, the composite index group
\[
\Gamma^\sharp=\Gamma\otimes_{\ZZ}T\otimes_{\ZZ}\Gamma,
\qquad
a\cdot(\alpha\otimes b\otimes\beta)\cdot c=a\alpha b\beta c,
\]
makes $T$ an ordinary Barnes $\Gamma^\sharp$-ring.  As proved in
\cref{prop:derived-gamma-comparison}, its raw left-operator algebra is the same
four-slot tensor rng $F(T,\Gamma)$, and its classical faithful left operator
rng is the nonunital regular image $\mathcal L_T$.  The new quotient between
them is precise:
our module notion imposes the additional middle-contraction identity, so that
the comparison factors as
\[
F(T,\Gamma)\twoheadrightarrow\cO(T,\Gamma)
\twoheadrightarrow\mathcal L_T\cong L(T,\Gamma^\sharp),
\qquad
\mathcal L_T\hookrightarrow\cR_T
  =\ZZ\,\id_T+\mathcal L_T.
\]
Thus the contribution is not the mere occurrence of four slots.  It is the
closed integral presentation that retains the decomposed composite index,
imposes the middle-module contraction, and exposes the kernel between that
universal quotient and the classical regular operator rng.

\begin{center}
\small
\begin{tabularx}{0.98\textwidth}{@{}>{\raggedright\arraybackslash}p{0.22\textwidth}>{\raggedright\arraybackslash}X>{\raggedright\arraybackslash}X@{}}
\toprule
 & Classical $\Gamma$-operator rings & Present ternary operator envelope \\ \midrule
Index group & $\Gamma^\sharp$ & decomposed $\Gamma\otimes T\otimes\Gamma$ \\
Raw generator & $a\otimes\xi$ & $a\otimes\alpha\otimes b\otimes\beta$ \\
Module law & outer associativity & outer and middle associativity \\
Resulting quotient & faithful nonunital regular operator rng $L$ & universal middle quotient $\cO$, then $\cO\twoheadrightarrow\mathcal L_T\cong L$ \\
Invariant studied & regular operator action & kernel of $\cO\to\mathcal L_T$ and its arithmetic; unitalization gives $\cR_T$ \\
\bottomrule
\end{tabularx}
\end{center}

\subsection{Ternary rings, semirings and semimodules}
Lister's ternary rings and the structural theory of associative triple systems provide the basic totally associative ternary identity \cite{Lister1971,Meyberg1972}.  Ternary semirings were introduced and developed through regularity, radicals, ideals and semifields by Dutta, Kar and collaborators \cite{DuttaKar2003,DuttaKar2004,DuttaKar2005,DuttaShumMandal2012}.  The literature also contains right and left ternary semimodules, annihilator constructions and singular ideals \cite{DuttaShumMandal2012,WaniPawar2016}.  These works are indispensable antecedents, but their ambient additive structures are semigroups and their principal tools are ideal-theoretic.  Our additive-group hypothesis is stronger and deliberate: it makes kernels, cokernels, tensor products, projective resolutions and Smith normal form available without subtractive replacements.  The operator defect is therefore not a renamed singular ideal; it is the kernel of a universal representation map between associative envelopes.

\subsection{Universal envelopes, standard embeddings and intrinsic triple cohomology}
Universal associative envelopes of triple systems usually start from the free associative algebra on the underlying triple space and force the ternary product to equal a polynomial in embedded generators \cite{Bremner2014,Elgendy2014}.  Standard embeddings likewise place a triple system inside a larger binary algebra.  Our universal object has a different representing problem: it universalizes the operators acting on all ternary modules, not an embedding of the elements of $T$ themselves.  Consequently the regular quotient may forget universal slot identities even when the ternary system embeds faithfully elsewhere.

There is also a general operadic enveloping-algebra principle.  For an algebra $V$ over a symmetric operad $\mathcal P$, the associative algebra $\mathsf U_{\mathcal P}(V)$ represents the corresponding module category; the construction and its PBW behavior are developed in \cite{LodayVallette2012,Khoroshkin2022}.  Accordingly, we do \emph{not} claim that passage from a multiary module law to some associative representing algebra is new in the abstract.  For the left-module notion fixed in \cref{def:module}, we prove that the representing rng has the closed integral presentation
\[
T\otimes\Gamma\otimes T\otimes\Gamma\big/J(T,\Gamma),
\]
and that the middle relations generate its contraction ideal.  The comparison
in \cref{prop:derived-gamma-comparison} identifies exactly how this quotient
differs from the operator ring of the derived Barnes $\Gamma^\sharp$-ring.
This explicit form retains the two $T$-slots and two $\Gamma$-slots needed for
Smith invariants, bad-prime fibers and the defect filtration; those arithmetic
data are not supplied by the abstract operadic existence theorem.  One may
regard $\cA(T,\Gamma)$ as the concrete enveloping algebra of the fixed-$\Gamma$
module theory, but the principal results concern the additional quotient
$\cA\twoheadrightarrow\cR_T$ and its invisible-operator arithmetic.

Carlsson's construction is a closer enveloping precedent than its cohomological
use alone suggests.  For an associative triple system $A$, Carlsson embeds $A$
in an associative algebra $U(A)$ and associates to a trimodule $M$ a
$U(A)$-bimodule $M_u$; Loos's standard embedding likewise places associative
triple systems in binary associative structures \cite{Carlsson1976,Carlsson1977,Loos1972}.
The specialization $\Gamma=\ZZ$ therefore has the following precise
comparison.
\begin{center}
\small
\begin{tabularx}{0.96\textwidth}{@{}>{\raggedright\arraybackslash}p{0.22\textwidth}>{\raggedright\arraybackslash}X>{\raggedright\arraybackslash}X@{}}
\toprule
 & Carlsson & Present construction \\ \midrule
Input & associative triple system $A$ & additive ternary $\Gamma$-ring $T$ \\
Module object & trimodule, with all module positions & one-sided last-slot module of \cref{def:module} \\
Associative object & embedding algebra $U(A)$ and the bimodule $M_u$ & representing algebra $\ZZ\ltimes\cO(T,\Gamma)$ \\
Invariant used here & intrinsic triple cohomology & kernel and extension of $\cA\twoheadrightarrow\cR_T$ \\
\bottomrule
\end{tabularx}
\end{center}
Because the module notions differ, we assert no canonical isomorphism between
$U(A)$ and $\ZZ\ltimes\cO(A,\ZZ)$.  The novelty claimed here is instead the
integral, $\Gamma$-indexed presentation for the specified one-sided theory and
the arithmetic of its regular kernel.  Likewise, the Hochschild class below is
not a replacement for intrinsic triple cohomology: it belongs to
\(0\to K\to\cA\to\cR_T\to0\) and measures whether universal operators split
from regularly visible operators.

\subsection{Conormal modules, filtered deformation and exact arithmetic}
For a ring epimorphism, the conormal bimodule $K/K^2$, its first relative
$\Tor$/$\Ext$ descriptions, the filtration $K^n/K^{n+1}$, and the tensor-algebra
map onto the associated graded are standard quotient-ring constructions
\cite{CartanEilenberg1956,Weibel1994}.  Hochschild cohomology classifies
additively split square-zero associative extensions
\cite{Hochschild1945,Witherspoon2019}, while Gerstenhaber's deformation theory
explains why filtered rings and their associated graded objects retain
successive obstruction data \cite{Gerstenhaber1964,Gerstenhaber1966}.  We use
these results as infrastructure.  The paper's specific content is that $K$ is
canonically produced by the ternary regular representation, together with the
relations among and computations of these invariants in the families below.

For finite-free tensor systems, flattenings and determinantal rank loci connect
the theory with the geometry of tensors \cite{Landsberg2012}.  Fitting ideals
and their base-change property encode the corresponding rank and support strata
\cite{Eisenbud1995,StacksFitting}, and Smith normal form records the integral
torsion exactly; polynomial-time algorithms for the latter were established by
Kannan and Bachem \cite{KannanBachem1979}.  These are standard tools.  The new
input is the two-matrix decomposition identifying which base-change obstruction
belongs to regular visibility and which torsion is already inherited from
middle relations.  In the sign family, classical balance and switching theory
for signed incidence matrices \cite{Zaslavsky1982} is applied to $B_\mu$; the
new output is the resulting exact $224+28+4$ census for this $256$-element
ensemble.

\subsection{Claimed contribution}
The paper's contribution is not the abstract existence of an associative
envelope, the classical operator ring of the derived $\Gamma^\sharp$-ring, or
the formal conormal, Hochschild, Fitting and matrix-factorization mechanisms.
To the authors' knowledge, the central new statement is the first realization
theorem prescribing independently the two coupled torsion cokernels attached
to a surjective integral tensor $E^{\otimes3}\twoheadrightarrow W$: two
complementary atomic families and an orthogonal-sum construction realize an
arbitrary ordered pair of finite abelian groups.  Realizing a group as the
cokernel of an arbitrary integer matrix does not imply this result, because
both presentation matrices must arise from the same trilinear tensor.  The
explicit two-step structure formula and the exact sign-tensor census are
supporting contributions; the standard base-change, determinantal and periodic
tools are applications that interpret these ternary operator quotients.

\section{Additive ternary \texorpdfstring{$\Gamma$}{Gamma}-rings and modules}

Throughout, $(T,+)$ and $(\Gamma,+)$ are abelian groups. The additive-group hypothesis is intentional: kernels, cokernels, tensor products over $\ZZ$, derived functors and Smith normal forms are then available without subtractivity hypotheses.  Here and below, a map described as additive in several variables is required to be \emph{separately additive in each displayed variable}; equivalently, it is $\ZZ$-multilinear.  It is not merely an additive homomorphism from the Cartesian product.  An \emph{rng} is associative but need not have an identity; rng homomorphisms are additive and multiplicative.  Rings and their modules are unital unless explicitly stated otherwise.

\begin{definition}\label{def:tgring}
An \emph{additive ternary $\Gamma$-ring} is an abelian group $T$ equipped with a separately additive map
\[
T\times\Gamma\times T\times\Gamma\times T\longrightarrow T,
\qquad
(x,\alpha,y,\beta,z)\longmapsto x\alpha y\beta z,
\]
satisfying
\begin{equation}\label{eq:ternary-assoc}
(x\alpha y\beta z)\gamma u\delta v
=x\alpha(y\beta z\gamma u)\delta v
=x\alpha y\beta(z\gamma u\delta v).
\end{equation}
\end{definition}

Because the operation is additive in every variable, inserting zero in any element or index slot gives zero.

\begin{definition}\label{def:tgring-morphism}
A morphism $(T,\Gamma)\to(T',\Gamma')$ of additive ternary indexed rings is a
pair $(f,g)$ of additive homomorphisms such that
\[
f(x\alpha y\beta z)
=f(x)\,g(\alpha)\,f(y)\,g(\beta)\,f(z).
\]
For a fixed index group $\Gamma$, morphisms mean pairs $(f,\id_\Gamma)$.
\end{definition}

\begin{definition}\label{def:module}
A \emph{left additive ternary $\Gamma$-module} over $T$ is an abelian group $M$ equipped with a separately additive action
\[
T\times\Gamma\times T\times\Gamma\times M\longrightarrow M,
\quad
(a,\alpha,b,\beta,m)\longmapsto a\alpha b\beta m,
\]
such that
\begin{equation}\label{eq:module-assoc}
(x\alpha y\beta z)\gamma u\delta m
=x\alpha(y\beta z\gamma u)\delta m
=x\alpha y\beta(z\gamma u\delta m).
\end{equation}
Morphisms are additive action-preserving maps. The category is denoted $T\text{-}\Gamma\mathrm{Mod}$.
\end{definition}

\begin{example}[Central-scalar systems]\label{ex:central}
Let $S$ be an associative rng and let $\lambda:\Gamma\to Z(S)$ be an additive homomorphism. Then
\[
x\alpha y\beta z=x\lambda(\alpha)y\lambda(\beta)z
\]
defines an additive ternary $\Gamma$-ring. Every left $S$-module becomes a ternary $\Gamma$-module by the analogous formula.
\end{example}

\begin{remark}\label{rem:singleton}
The one-element additive group $\Gamma=0$ forces the indexed product to be zero. Nonzero ordinary ternary products are recovered, for example, by taking $\Gamma=\ZZ$ and setting $x\,n\,y\,m\,z=nm[x,y,z]$. This avoids the invalid use of a unique zero index as if it were a multiplicative unit.
\end{remark}

\section{The universal operator rng}

Set
\[
F(T,\Gamma)=T\otimes_{\ZZ}\Gamma\otimes_{\ZZ}T\otimes_{\ZZ}\Gamma.
\]
A pure tensor is written $a\otimes\alpha\otimes b\otimes\beta$.

\begin{lemma}\label{lem:Fmult}
There is a unique bilinear multiplication on $F(T,\Gamma)$ satisfying
\begin{equation}\label{eq:Fmult}
(a\otimes\alpha\otimes b\otimes\beta)
(c\otimes\gamma\otimes d\otimes\delta)
=(a\alpha b\beta c)\otimes\gamma\otimes d\otimes\delta.
\end{equation}
It is associative.
\end{lemma}

\begin{proof}
Consider first the eight-variable map on pure slots whose value is the right-hand side of \eqref{eq:Fmult}.  Separate additivity of the ternary product and of the canonical pure-tensor map makes this map separately additive in all eight variables.  Applying the universal property of the tensor product in the first four variables gives, for each fixed $(c,\gamma,d,\delta)$, an additive map $F(T,\Gamma)\to F(T,\Gamma)$.  Its dependence on the remaining four variables is again separately additive, so a second application of the universal property produces a unique bilinear map
\[
F(T,\Gamma)\times F(T,\Gamma)\longrightarrow F(T,\Gamma)
\]
with the displayed pure-tensor formula.  Thus the proposed multiplication is well defined; this argument also records the balancing in each tensor slot.

On three pure tensors the first tensor factors obtained from the two parenthesizations are
\[
(a\alpha b\beta c)\gamma d\delta e
\quad\text{and}\quad
a\alpha b\beta(c\gamma d\delta e),
\]
which agree by \eqref{eq:ternary-assoc}.  Since pure tensors generate $F(T,\Gamma)$ additively and the multiplication is bilinear, the equality extends to arbitrary three elements of $F(T,\Gamma)$.
\end{proof}

The middle parenthesization in \eqref{eq:module-assoc} imposes an additional family of relations.

\begin{definition}\label{def:O}
Let $J(T,\Gamma)$ be the two-sided ideal generated by
\begin{equation}\label{eq:middle}
(a\alpha b\beta c)\otimes\gamma\otimes d\otimes\delta
-a\otimes\alpha\otimes(b\beta c\gamma d)\otimes\delta.
\end{equation}
The \emph{universal operator rng} is
\[
\cO(T,\Gamma)=F(T,\Gamma)/J(T,\Gamma).
\]
The class of a pure tensor is denoted $[a,\alpha,b,\beta]$.
\end{definition}

\begin{proposition}[Structure of the relation ideal]\label{prop:relation-ideal}
Let $N(T,\Gamma)$ be the additive subgroup of $F(T,\Gamma)$ generated by the elements in \eqref{eq:middle}. Then
\[
N(T,\Gamma)F(T,\Gamma)=0,
\qquad
F(T,\Gamma)N(T,\Gamma)\subseteq N(T,\Gamma).
\]
Consequently $N(T,\Gamma)$ is already a two-sided ideal and $J(T,\Gamma)=N(T,\Gamma)$.
\end{proposition}

\begin{proof}
Write $r(a,\alpha,b,\beta,c,\gamma,d,\delta)$ for the generator in \eqref{eq:middle}.  Multiplication on the right by a pure operator $e\otimes\varepsilon\otimes f\otimes\zeta$ gives zero, because its two first tensor factors are
\[
(a\alpha b\beta c)\gamma d\delta e
\quad\text{and}\quad
a\alpha(b\beta c\gamma d)\delta e,
\]
which agree by the middle associativity identity.  Multiplication on the left gives
\[
(e\varepsilon f\zeta a)\alpha b\beta c\otimes\gamma\otimes d\otimes\delta
-(e\varepsilon f\zeta a)\otimes\alpha\otimes(b\beta c\gamma d)\otimes\delta,
\]
which is $r(e\varepsilon f\zeta a,\alpha,b,\beta,c,\gamma,d,\delta)$; here the right associativity identity was used in the first term.  Pure tensors generate $F(T,\Gamma)$ and multiplication is bilinear, so both conclusions extend to the subgroup generated by all $r$.  Thus $N$ is a two-sided ideal containing the defining generators, and minimality of $J$ gives $J=N$.
\end{proof}

\begin{definition}\label{def:realization}
An \emph{operator realization} of $(T,\Gamma)$ in an associative rng $B$ is a map
\[
\lambda:T\times\Gamma\times T\times\Gamma\to B
\]
additive in all variables and satisfying
\begin{align}
\lambda(a,\alpha,b,\beta)\lambda(c,\gamma,d,\delta)
&=\lambda(a\alpha b\beta c,\gamma,d,\delta),\label{eq:real1}\\
\lambda(a\alpha b\beta c,\gamma,d,\delta)
&=\lambda(a,\alpha,b\beta c\gamma d,\delta).\label{eq:real2}
\end{align}
\end{definition}

\begin{proposition}[Universal property]\label{thm:universal}
The map
\[
\iota:T\times\Gamma\times T\times\Gamma\to\cO(T,\Gamma),
\qquad
(a,\alpha,b,\beta)\mapsto[a,\alpha,b,\beta],
\]
is an operator realization. For every operator realization $\lambda$ in an associative rng $B$, there is a unique rng homomorphism
\[
\overline\lambda:\cO(T,\Gamma)\to B
\]
satisfying $\overline\lambda\circ\iota=\lambda$.
\end{proposition}

\begin{proof}
The defining pure-class map $\iota$ is separately additive.  The multiplication
formula in \cref{lem:Fmult} gives \eqref{eq:real1} for $\iota$, while the
vanishing of \eqref{eq:middle} in the quotient gives \eqref{eq:real2}.  Thus
$\iota$ is an operator realization.

Now let $\lambda$ be an operator realization in $B$.  Four successive tensor
universal properties give a unique additive map
\[
\widetilde\lambda:F(T,\Gamma)\longrightarrow B,
\qquad
\widetilde\lambda(a\otimes\alpha\otimes b\otimes\beta)
=\lambda(a,\alpha,b,\beta).
\]
For pure $x,y\in F(T,\Gamma)$, identity \eqref{eq:real1} and
\eqref{eq:Fmult} give
$\widetilde\lambda(xy)=\widetilde\lambda(x)\widetilde\lambda(y)$.
For fixed $y$, both sides are additive in $x$, so tensor extensionality extends
the equality to arbitrary $x$; repeating the argument in $y$ proves that
$\widetilde\lambda$ is an rng homomorphism on all of $F(T,\Gamma)$.

Identity \eqref{eq:real2} says exactly that $\widetilde\lambda$ sends every
generator \eqref{eq:middle} to zero.  Additivity and
\cref{prop:relation-ideal} therefore give
$N(T,\Gamma)=J(T,\Gamma)\subseteq\ker\widetilde\lambda$.  Hence
$\widetilde\lambda(x)=\widetilde\lambda(y)$ whenever $x-y\in J(T,\Gamma)$,
so it induces an rng homomorphism
$\overline\lambda:\cO(T,\Gamma)\to B$ satisfying
$\overline\lambda\circ\iota=\lambda$.

Finally, if $f:\cO(T,\Gamma)\to B$ has the same property, then
$f\circ\pi$ and $\widetilde\lambda$ are additive maps on $F(T,\Gamma)$ that
agree on every pure tensor, where $\pi:F(T,\Gamma)\to\cO(T,\Gamma)$ is the
quotient map.  Tensor extensionality gives
$f\circ\pi=\widetilde\lambda$; surjectivity of $\pi$ then yields
$f=\overline\lambda$.
\end{proof}

\begin{proposition}[Comparison with the derived Barnes $\Gamma$-ring]
\label{prop:derived-gamma-comparison}
Put
\[
\Gamma^\sharp=\Gamma\otimes_{\ZZ}T\otimes_{\ZZ}\Gamma
\]
and define, first on pure tensors,
\begin{equation}\label{eq:derived-gamma-product}
a\cdot(\alpha\otimes b\otimes\beta)\cdot c
=a\alpha b\beta c.
\end{equation}
Then $(T,\Gamma^\sharp)$ is a Barnes $\Gamma$-ring.  Under tensor
reassociation, its raw left-operator rng is canonically isomorphic to
$F(T,\Gamma)$.  Let
\[
\lambda_T:F(T,\Gamma)\longrightarrow\End_{\ZZ}(T),
\qquad
\lambda_T(a\otimes\alpha\otimes b\otimes\beta)(x)=a\alpha b\beta x,
\]
and put $\mathcal L_T=\im\lambda_T$.  If $L(T,\Gamma^\sharp)$ denotes the
classical left operator ring obtained by identifying finite sums that induce
the same endomorphism of the regular module, then there are canonical rng maps
\begin{equation}\label{eq:derived-gamma-factorization}
F(T,\Gamma)\twoheadrightarrow\cO(T,\Gamma)
\twoheadrightarrow\mathcal L_T\xrightarrow{\ \cong\ }L(T,\Gamma^\sharp).
\end{equation}
The kernel of the second map is the intrinsic regular-action kernel, denoted
$\Def(T,\Gamma)$ in \cref{def:defects}.  Its unital enlargement satisfies
$\cR_T=\ZZ\,\id_T+\mathcal L_T$ inside $\End_{\ZZ}(T)$.

Moreover, an ordinary left module over $(T,\Gamma^\sharp)$ requires the outer
identity
\[
(a\cdot\xi\cdot c)\cdot\eta\cdot m
=a\cdot\xi\cdot(c\cdot\eta\cdot m).
\]
Under \eqref{eq:derived-gamma-product}, the modules represented by
$\cO(T,\Gamma)$ form the full subcategory in which the additional middle
identity
\[
(a\alpha b\beta c)\gamma d\delta m
=a\alpha(b\beta c\gamma d)\delta m
\]
also holds.  Thus the middle-relation ideal $J(T,\Gamma)$ is exactly the
additional contraction imposed by the module notion of \cref{def:module}.
\end{proposition}

\begin{proof}
Separate additivity of the ternary product extends
\eqref{eq:derived-gamma-product} uniquely from pure tensors to
$\Gamma^\sharp$.  For pure
$\xi=\alpha\otimes b\otimes\beta$ and
$\eta=\gamma\otimes d\otimes\delta$, the Barnes identity is
\[
(a\cdot\xi\cdot c)\cdot\eta\cdot e
=(a\alpha b\beta c)\gamma d\delta e
=a\alpha b\beta(c\gamma d\delta e)
=a\cdot\xi\cdot(c\cdot\eta\cdot e),
\]
and multilinearity gives the general case.

The tensor reassociation
\[
T\otimes\Gamma^\sharp\longrightarrow F(T,\Gamma),
\qquad
a\otimes(\alpha\otimes b\otimes\beta)
\longmapsto a\otimes\alpha\otimes b\otimes\beta
\]
is an additive isomorphism.  The classical pair-operator product sends
$(a\otimes\xi)(c\otimes\eta)$ to
$(a\cdot\xi\cdot c)\otimes\eta$; hence the displayed isomorphism is
multiplicative by \eqref{eq:Fmult}.  Its regular action on $x\in T$ is
$x\mapsto a\alpha b\beta x$, exactly the map defining $\mathcal L_T$.
Quotienting by equality of regular actions therefore gives
$L(T,\Gamma^\sharp)\cong\mathcal L_T$.

Every generator \eqref{eq:middle} acts trivially on $T$ by the middle equality
in \eqref{eq:ternary-assoc}; hence the regular map factors through $\cO$, and
its kernel there is the intrinsic regular-action kernel.  Adjoining integer
scalar operators gives $\cR_T=\ZZ\,\id_T+\mathcal L_T$.  Finally, the ordinary
derived-$\Gamma^\sharp$ module identity is the first--third equality in
\eqref{eq:module-assoc}; imposing its first--middle equality kills precisely
the generators \eqref{eq:middle}.  The universal property
\cref{thm:universal} proves the asserted full-subcategory statement.
\end{proof}

\begin{remark}[Relation to other envelopes]\label{rem:other-envelopes}
Universal associative envelopes of nonassociative triple systems generally begin with a free associative algebra on the underlying triple system and impose relations that express the ternary product as a polynomial in embedded generators; see \cite{Bremner2014,Elgendy2014}. The object $\cO(T,\Gamma)$ is different. Its generators are two-element, two-index operator slots, and its universal property concerns representations on ternary modules. It is therefore closer in spirit to an operator ring or standard embedding, while retaining the $\Gamma$-indexed middle relation.
\end{remark}

\begin{proposition}\label{prop:module-rep}
Every ternary $\Gamma$-module $M$ determines a unique rng homomorphism
\[
\rho_M:\cO(T,\Gamma)\to\End_{\ZZ}(M),
\qquad
\rho_M([a,\alpha,b,\beta])(m)=a\alpha b\beta m.
\]
This construction is functorial: if $f:M\to N$ is a ternary
$\Gamma$-module morphism, then
\[
f\,\rho_M(u)=\rho_N(u)\,f
\qquad (u\in\cO(T,\Gamma)).
\]
\end{proposition}

\begin{proof}
Give $\End_{\ZZ}(M)$ its usual multiplication by composition, so that
$(FG)(m)=F(G(m))$.  The separately additive map
\[
(a,\alpha,b,\beta)\longmapsto
\bigl(m\mapsto a\alpha b\beta m\bigr)
\]
satisfies \eqref{eq:real1}: after evaluation at $m$, the required identity is
the equality of the first and third expressions in \eqref{eq:module-assoc}.
It satisfies \eqref{eq:real2} by the equality of the first and second
expressions there.  Hence it is an operator realization, and
\cref{thm:universal} gives the unique rng homomorphism $\rho_M$ with the
displayed action formula.

For functoriality, let $f:M\to N$ be a ternary-module morphism.  On a pure
operator $[a,\alpha,b,\beta]$, action preservation gives
\[
f\rho_M([a,\alpha,b,\beta])(m)
=f(a\alpha b\beta m)
=a\alpha b\beta f(m)
=\rho_N([a,\alpha,b,\beta])f(m).
\]
For fixed $f$ and $m$, both sides are additive in $u$.  Pure operator classes
generate $\cO(T,\Gamma)$ additively, because they generate the fourfold tensor
group before quotienting.  The intertwining identity therefore holds for
every $u\in\cO(T,\Gamma)$.
\end{proof}

\section{Dorroh extension and categorical representation}

For an associative rng $R$, its Dorroh extension $\ZZ\ltimes R$ is the abelian group $\ZZ\oplus R$ with
\begin{equation}\label{eq:dorroh}
(n,r)(m,s)=(nm,ns+mr+rs).
\end{equation}
It is a unital ring with identity $(1,0)$ \cite{Dorroh1932}.

\begin{definition}\label{def:A}
Set
\[
\cA(T,\Gamma)=\ZZ\ltimes\cO(T,\Gamma).
\]
\end{definition}

\begin{proposition}[Operator equivalence]\label{thm:equiv}
The category $T\text{-}\Gamma\mathrm{Mod}$ is canonically equivalent to the category of unital left $\cA(T,\Gamma)$-modules.
\end{proposition}

\begin{proof}
Write $\iota(u)=(0,u)$ and $\eta(n)=(n,0)$.  These maps satisfy
\begin{equation}\label{eq:dorroh-generators}
\eta(n)\eta(m)=\eta(nm),\qquad
\eta(n)\iota(u)=\iota(nu),\qquad
\iota(u)\eta(n)=\iota(nu),\qquad
\iota(u)\iota(v)=\iota(uv),
\end{equation}
and every element of $\cA$ has the unique decomposition
$\eta(n)+\iota(u)$.

Let $M$ be a ternary module.  By \cref{prop:module-rep}, its pure actions
extend uniquely to a rng homomorphism
$\rho_M:\cO\to\End_{\ZZ}(M)$.  Define
\[
\widetilde\rho_M(n,u)(m)=nm+\rho_M(u)(m).
\]
For $x=(n,u)$ and $y=(m,v)$, additivity of the endomorphisms and
multiplicativity of $\rho_M$ give
\begin{align*}
\widetilde\rho_M(x)\widetilde\rho_M(y)(z)
 &=nmz+n\rho_M(v)z+m\rho_M(u)z+\rho_M(uv)z\\
 &=\widetilde\rho_M\bigl(nm,nv+mu+uv\bigr)(z).
\end{align*}
Moreover $\widetilde\rho_M(1,0)=\id_M$.  Hence
$\widetilde\rho_M$ is a unital ring homomorphism and makes $M$ a unital
$\cA$-module.

Conversely, let $M$ be a unital left $\cA$-module, with representation
$\sigma:\cA\to\End_{\ZZ}(M)$.  Restricting $\sigma$ along $\iota$ gives a rng
homomorphism $\rho_\sigma:\cO\to\End_{\ZZ}(M)$.  Define
\[
a\alpha b\beta m
=\rho_\sigma([a,\alpha,b,\beta])(m)
=\sigma(0,[a,\alpha,b,\beta])(m).
\]
Separate additivity follows from the fourfold tensor map.  Multiplication of
pure classes gives
\[
[a,\alpha,b,\beta]\,[c,\gamma,d,\delta]
=[a\alpha b\beta c,\gamma,d,\delta],
\]
while the defining middle relation gives
\[
[a\alpha b\beta c,\gamma,d,\delta]
=[a,\alpha,b\beta c\gamma d,\delta].
\]
After applying $\rho_\sigma$ and evaluating at $m$, these are precisely the
right and middle equalities in \eqref{eq:module-assoc}; thus the recovered
operation is a ternary-module action.

The object constructions are inverse.  Starting with a ternary action,
restriction of $\widetilde\rho_M$ to $\iota(\cO)$ is $\rho_M$, so the pure
actions are unchanged.  Starting with $\sigma$, the reconstruction has
operator part $\rho_\sigma$, and the re-extended representation satisfies
\[
\widetilde\rho_{\sigma}(n,u)
=n\,\id_M+\sigma(0,u)
=\sigma(n,0)+\sigma(0,u)
=\sigma(n,u).
\]

It remains to check morphisms.  An additive map $f:M\to N$ between two
ternary modules preserves the pure actions if and only if it intertwines
$\rho_M$ and $\rho_N$ on the pure operator classes.  Those classes generate
$\cO$ additively, so this is equivalent to intertwining every $u\in\cO$.
Every additive map commutes with the integer-scalar operators; using the
unique decomposition $(n,u)=(n,0)+(0,u)$, the latter condition is equivalent
to
\[
f\,\widetilde\rho_M(x)=\widetilde\rho_N(x)f
\qquad(x\in\cA).
\]
Thus the two object maps induce mutually inverse bijections on every Hom-set,
and hence define a canonical equivalence of categories.
\end{proof}

\begin{corollary}\label{cor:abelian}
The category $T\text{-}\Gamma\mathrm{Mod}$ is abelian and has enough projectives and injectives.
\end{corollary}

\subsection{Scalar extension of the envelope}

Let $S$ be a commutative unital $\ZZ$-algebra.  Put
\[
T_S=T\otimes_{\ZZ}S,
\qquad
\Gamma_S=\Gamma\otimes_{\ZZ}S,
\]
and extend the ternary product $S$-multilinearly; on pure tensors this means
\[
(x\otimes s_1)(\alpha\otimes s_2)(y\otimes s_3)
(\beta\otimes s_4)(z\otimes s_5)
=(x\alpha y\beta z)\otimes s_1s_2s_3s_4s_5.
\]
Let
\[
F_S=T_S\otimes_S\Gamma_S\otimes_S T_S\otimes_S\Gamma_S
\]
with the raw multiplication given by the scalar-extended version of
\eqref{eq:Fmult}, let $J_S$ be the $S$-submodule generated by the same middle
relations, and put
\[
\cO_S(T_S,\Gamma_S)=F_S/J_S,
\qquad
\cA_S(T_S,\Gamma_S)=S\ltimes\cO_S(T_S,\Gamma_S).
\]
The proof below includes the fact that $J_S$ is a two-sided ideal, so these
quotients carry the asserted multiplicative structures.

\begin{proposition}[Envelope base change]\label{thm:envelope-basechange}
For every commutative unital $\ZZ$-algebra $S$, scalar extension induces natural isomorphisms
\begin{align*}
\cO(T,\Gamma)\otimes_{\ZZ}S
&\cong \cO_S(T_S,\Gamma_S),\\
\cA(T,\Gamma)\otimes_{\ZZ}S
&\cong \cA_S(T_S,\Gamma_S)
\end{align*}
of $S$-rngs and unital $S$-algebras, respectively.
\end{proposition}

\begin{proof}
Write $F=F(T,\Gamma)$, $J=J(T,\Gamma)$ and $\cO=F/J$, and let
$j:J\hookrightarrow F$ and $q:F\twoheadrightarrow\cO$ be the inclusion and
quotient maps.  Tensor-product reassociation gives a canonical $S$-linear
isomorphism
\[
\Phi:F\otimes_{\ZZ}S\xrightarrow{\ \sim\ }F_S.
\]
For example, $\Phi((a\otimes\alpha\otimes b\otimes\beta)\otimes s)$
may be written as
$(a\otimes s)\otimes_S(\alpha\otimes1)\otimes_S
(b\otimes1)\otimes_S(\beta\otimes1)$.

Independently of $\Phi$, give $F\otimes_{\ZZ}S$ the bilinear multiplication
\[
(u\otimes s)(v\otimes t)=uv\otimes st
\qquad(u,v\in F,\ s,t\in S).
\]
The tensor universal property defines this multiplication on all elements,
and associativity follows from associativity in $F$ and $S$.  Formula
\eqref{eq:Fmult} on pure tensors shows that $\Phi$ is multiplicative.

Tensoring the exact presentation $J\xrightarrow{j}F\xrightarrow{q}\cO\to0$
with $S$ gives the right-exact sequence
\[
J\otimes_{\ZZ}S\xrightarrow{j\otimes1}
F\otimes_{\ZZ}S\xrightarrow{q\otimes1}
\cO\otimes_{\ZZ}S\longrightarrow0.
\]
In particular,
\begin{equation}\label{eq:basechange-range}
\ker(q\otimes1)=\operatorname{im}(j\otimes1).
\end{equation}
This is an image statement: $j\otimes1$ need not be injective because no
flatness assumption is made on $S$.

We claim
\begin{equation}\label{eq:basechange-relations}
\Phi\bigl(\operatorname{im}(j\otimes1)\bigr)=J_S.
\end{equation}
To verify this without suppressing the tensor balancing, write
\begin{align*}
\mathfrak r_R(a,\alpha,b,\beta,c,\gamma,d,\delta)
={}&(a\alpha b\beta c)\otimes_R\gamma\otimes_R d\otimes_R\delta\\
&-a\otimes_R\alpha\otimes_R(b\beta c\gamma d)\otimes_R\delta
\end{align*}
for a displayed middle relation over $R$, where $R=\ZZ$ or $S$.  The
pure-generator calculation is
\begin{multline}\label{eq:pure-relation-basechange}
\mathfrak r_S(a\otimes s_1,\alpha\otimes s_2,
b\otimes s_3,\beta\otimes s_4,c\otimes s_5,\gamma\otimes s_6,
d\otimes s_7,\delta\otimes s_8)\\
=\Phi\!\left(
\mathfrak r_{\ZZ}(a,\alpha,b,\beta,c,\gamma,d,\delta)
\otimes s_1s_2s_3s_4s_5s_6s_7s_8\right).
\end{multline}
Indeed, the two sides expand to the same difference of pure fourfold
tensors by the displayed formula for the extended ternary product and
$S$-balancing.  Taking $(s_1,\ldots,s_8)=(s,1,\ldots,1)$ proves that the
transport of every integral generator tensored with $s$ lies in $J_S$.
Conversely, expand each entry of an arbitrary middle relation over $S$ as a
sum of pure tensors.  Eight-slot additivity reduces every summand to the
left-hand side of \eqref{eq:pure-relation-basechange}, which belongs to the
transported image by its right-hand side.  Since both sides of
\eqref{eq:basechange-relations} are $S$-submodules, the two generator
inclusions prove the claimed equality.

Because $J$ is a two-sided ideal of $F$, the image
$\operatorname{im}(j\otimes1)$ is a two-sided ideal of
$F\otimes_{\ZZ}S$: multiplication by a pure tensor sends
$j(x)\otimes s$ to $(uj(x))\otimes ts$ or $(j(x)u)\otimes st$, and
bilinearity gives the general case.  Since $\Phi$ is multiplicative,
\eqref{eq:basechange-relations} shows in particular that $J_S$ is a two-sided
ideal of $F_S$.  Equations \eqref{eq:basechange-range} and
\eqref{eq:basechange-relations} therefore induce a multiplicative $S$-linear
isomorphism
\[
\overline\Phi:\cO\otimes_{\ZZ}S\xrightarrow{\sim}F_S/J_S=\cO_S.
\]
This proves the first isomorphism.

For the unital envelope, use the explicit additive isomorphism
\[
(\ZZ\ltimes\cO)\otimes_{\ZZ}S
\xrightarrow{\ \sim\ } S\oplus(\cO\otimes_{\ZZ}S),
\qquad (n,u)\otimes s\longmapsto(ns,u\otimes s).
\]
Under this map the product becomes
\[
(s,u)(t,v)=(st,\,sv+tu+uv),
\]
which is precisely the $S$-Dorroh multiplication.  Composing the second
coordinate with $\overline\Phi$ proves the asserted algebra isomorphism.
All maps used above are defined by tensor reassociation and quotient universal
properties, so the two isomorphisms commute with morphisms of ternary
$\Gamma$-rings and with homomorphisms of scalar algebras.  Hence they are
natural.
\end{proof}

The regular action over $S$ induces canonical maps
\[
\rho^S_{T_S}:\cO_S(T_S,\Gamma_S)\longrightarrow\End_S(T_S),
\qquad
\widehat\rho^S_{T_S}:\cA_S(T_S,\Gamma_S)\longrightarrow\End_S(T_S),
\]
where
$\rho^S_{T_S}([a,\alpha,b,\beta])(x)=a\alpha b\beta x$ and the second map sends $(s,u)$ to
$s\id_{T_S}+\rho^S_{T_S}(u)$.  We use the notation
\[
\begin{aligned}
\Def_S(T_S,\Gamma_S)&=\Ker\rho^S_{T_S},
&\quad \UDef_S(T_S,\Gamma_S)&=\Ker\widehat\rho^S_{T_S},\\
\mathcal L^S_{T_S}&=\im\rho^S_{T_S},
&\quad \mathcal R^S_{T_S}&=\im\widehat\rho^S_{T_S}.
\end{aligned}
\]
When $\Gamma=\ZZ$, so that $\Gamma_S\cong S$, we abbreviate
$\cO_S(T_S,\Gamma_S)$, $\cA_S(T_S,\Gamma_S)$ and
$\Def_S(T_S,\Gamma_S)$ to $\cO_S(T_S)$, $\cA_S(T_S)$ and
$\Def_S(T_S)$, respectively, and similarly for $\UDef_S$.

\begin{remark}\label{rem:envelope-versus-defect-basechange}
Thus the universal envelope itself has no scalar-extension obstruction.  Any failure of the regular defect to commute with base change comes from taking the kernel of the regular representation.  The exact correction term is computed for two-step tensor systems in \cref{thm:basechange}.
\end{remark}

\section{Regular images and operator-faithfulness defects}

The regular ternary module gives
\[
\rho_T:\cO(T,\Gamma)\to\End_{\ZZ}(T),
\qquad
\rho_T([a,\alpha,b,\beta])(x)=a\alpha b\beta x,
\]
and the unital extension
\[
\widehat\rho_T:\cA(T,\Gamma)\to\End_{\ZZ}(T),
\qquad
\widehat\rho_T(n,u)=n\id_T+\rho_T(u).
\]

\medskip
\begin{definition}\label{def:defects}
Define
\[
\Def(T,\Gamma)=\Ker\rho_T,
\qquad
\UDef(T,\Gamma)=\Ker\widehat\rho_T,
\qquad
\mathcal L_T=\im\rho_T,
\qquad
\cR_T=\im\widehat\rho_T.
\]
We call $\Def$ the \emph{intrinsic operator-faithfulness defect}, $\UDef$ the
\emph{unital regular defect}, $\mathcal L_T$ the \emph{raw regular operator
rng}, and $\cR_T$ the \emph{unital regular operator ring}.  Inside
$\End_{\ZZ}(T)$ one has
\[
\cR_T=\ZZ\,\id_T+\mathcal L_T;
\]
the sum need not be direct.
\end{definition}

\begin{proposition}\label{prop:def-seq}
Let
\[
I_T=\{n\in\ZZ:n\id_T\in\rho_T(\cO(T,\Gamma))\}.
\]
There is a canonical exact sequence
\begin{equation}\label{eq:def-seq}
0\longrightarrow\Def(T,\Gamma)
\longrightarrow\UDef(T,\Gamma)
\xrightarrow{\nu} I_T
\longrightarrow0,
\end{equation}
where $\nu(n,u)=n$.
\end{proposition}

\begin{proof}
If $(n,u)\in\UDef$, then $\rho_T(u)=-n\id_T$, so $n\in I_T$. The kernel is naturally $\Def$. Conversely, if $n\id_T=\rho_T(v)$, then $(n,-v)\in\UDef$.
\end{proof}

\begin{definition}
A ternary module is \emph{regularly compatible} if the corresponding $\cA(T,\Gamma)$-module is annihilated by $\UDef(T,\Gamma)$.
\end{definition}

\begin{proposition}\label{prop:regular-subcat}
Regularly compatible ternary modules form a full subcategory equivalent to $\cR_T\text{-}\mathrm{Mod}$.
\end{proposition}

\begin{proof}
The first isomorphism theorem gives $\cR_T\cong\cA/\UDef$. Modules over a quotient are precisely modules annihilated by its kernel.
\end{proof}

\section{The conormal defect and first relative homology}

For this section write
\[
\cA=\cA(T,\Gamma),\qquad
\cR=\cR_T,\qquad
K=\UDef(T,\Gamma).
\]
Thus
\begin{equation}\label{eq:quotient}
0\longrightarrow K\longrightarrow\cA\longrightarrow\cR\longrightarrow0
\end{equation}
is exact.

\begin{definition}\label{def:conormal}
The \emph{conormal operator defect} is the $\cR$-bimodule
\[
\cC(T,\Gamma)=K/K^2.
\]
\end{definition}

The $\cR$-bimodule structure is well defined because changing a lift of an element of $\cR$ changes its action on $K$ by a product in $K^2$.

\begin{proposition}[Conormal formula]\label{thm:conormal}
Let $M$ be a left $\cR$-module, viewed as a left $\cA$-module through the quotient. Then there are natural isomorphisms
\begin{align}
\Tor^{\cA}_1(\cR,M)
&\cong \cC(T,\Gamma)\otimes_{\cR}M,\label{eq:conormal-tor}\\
\Ext^1_{\cA}(\cR,M)
&\cong \Hom_{\cR}(\cC(T,\Gamma),M).\label{eq:conormal-ext}
\end{align}
Here $\Tor^{\cA}(\cR,M)$ uses $\cR$ as a right $\cA$-module and $M$ as a
left $\cA$-module, whereas $\Ext_{\cA}(\cR,M)$ uses $\cR$ and $M$ as left
$\cA$-modules.  On the right-hand sides, $\cC$ uses its right $\cR$-action in
the tensor product and its left $\cR$-action in the Hom group.
\end{proposition}

\begin{proof}
Applying $-\otimes_{\cA}M$ to \eqref{eq:quotient}, considered as right modules, gives
\[
0\to\Tor^{\cA}_1(\cR,M)\to K\otimes_{\cA}M\to\cA\otimes_{\cA}M.
\]
The last map sends $k\otimes m$ to $km=0$, so $\Tor_1\cong K\otimes_{\cA}M$. Balancedness kills $K^2$: for $k,k'\in K$,
\[
kk'\otimes m=k\otimes k'm=0.
\]
The remaining action factors through $\cR$, giving
\[
K\otimes_{\cA}M\cong(K/K^2)\otimes_{\cR}M.
\]

Applying $\Hom_{\cA}(-,M)$ to \eqref{eq:quotient} gives
\[
\Hom_{\cA}(\cA,M)\longrightarrow\Hom_{\cA}(K,M)\longrightarrow\Ext^1_{\cA}(\cR,M)\longrightarrow0.
\]
Every map $\cA\to M$ restricts to zero on $K$, because $K$ annihilates $M$. Hence $\Ext^1\cong\Hom_{\cA}(K,M)$. Any such map annihilates $K^2$ and is exactly an $\cR$-linear map $K/K^2\to M$.
\end{proof}

\begin{corollary}\label{cor:C-tor}
There is a canonical isomorphism of $\cR$-bimodules
\[
\cC(T,\Gamma)\cong\Tor^{\cA}_1(\cR,\cR).
\]
\end{corollary}

\begin{corollary}\label{cor:flat-projective}
If $\cC(T,\Gamma)\ne0$, then $\cR_T$ is not flat as a right $\cA(T,\Gamma)$-module.\par
If $\Hom_{\cR_T}(\cC(T,\Gamma),M)\ne0$ for some $M$, then $\cR_T$ is not projective as a left $\cA(T,\Gamma)$-module.
\end{corollary}

\section{The graded defect tower}

The first conormal layer need not capture the full ideal $K$. Its powers produce higher-order operator identities.

\begin{definition}\label{def:tower}
For $n\ge1$, define
\[
\cC_n(T,\Gamma)=K^n/K^{n+1}.
\]
The \emph{graded defect algebra} is
\[
\gr_K(\cA)=\cR\oplus\bigoplus_{n\ge1}\cC_n(T,\Gamma).
\]
\end{definition}

\begin{proposition}\label{prop:tower}
Each $\cC_n(T,\Gamma)$ is an $\cR$-bimodule. Multiplication induces canonical surjections
\[
\cC(T,\Gamma)^{\otimes_{\cR}n}\twoheadrightarrow\cC_n(T,\Gamma).
\]
Consequently, $\gr_K(\cA)$ is generated as a graded $\cR$-algebra by its degree-one part $\cC(T,\Gamma)$.
\end{proposition}

\begin{proof}
The $\cR$-action is defined by lifting elements to $\cA$; changing a lift by an element of $K$ changes the product by an element of $K^{n+1}$. Multiplication $K^{\otimes n}\to K^n$ is surjective by definition and factors through both $K/K^2$ and the balancing relations over $\cR$.
\end{proof}

\begin{definition}
If $K=0$, its \emph{operator-defect depth} is $0$.  If $K^{h+1}=0$ and $K^h\ne0$ for some $h\ge1$, its depth is $h$.  If $K$ is not nilpotent, its depth is infinite.
\end{definition}

The square-zero case has depth one. The rank-one family in \cref{sec:rank-one} has infinite depth but every graded layer is the same finite cyclic group.

\section{Higher relation modules and conormal rigidity}\label{sec:relations}

The surjections in \cref{prop:tower} have kernels that record relations among first-order invisible operators.

\begin{definition}\label{def:relationmodules}
For $n\ge2$, let
\[
\tau_n:\cC(T,\Gamma)^{\otimes_{\cR_T}n}
\longrightarrow \cC_n(T,\Gamma)
\]
be the multiplication map and define the $n$th \emph{higher operator-relation module}
\[
\cS_n(T,\Gamma)=\Ker\tau_n.
\]
The defect tower is called \emph{tensor-free} if every $\tau_n$ is an isomorphism.
\end{definition}

\begin{proposition}[Canonical presentation of the tower]\label{prop:relationpresentation}
The graded defect algebra has a canonical presentation
\[
\gr_K(\cA)
\cong
T_{\cR}(\cC)\big/\cJ_{\mathrm{gr}},
\]
where $T_{\cR}(\cC)=\cR\oplus\cC\oplus\cC^{\otimes2}\oplus\cdots$ is the tensor algebra and the homogeneous degree-$n$ component of $\cJ_{\mathrm{gr}}$ is $\cS_n$.  In particular, $\cS_2$ is the complete quadratic relation module of the regular defect filtration.
\end{proposition}

\begin{proof}
The direct sum of the multiplication maps $\tau_n$, together with the identity in degrees zero and one, is a surjective graded algebra homomorphism from the tensor algebra to $\gr_K(\cA)$.  Its homogeneous kernel in degree $n$ is precisely $\Ker\tau_n$.
\end{proof}

\begin{proposition}[Conormal rigidity for finite depth]\label{thm:conormal-rigidity}
If the unital regular defect $K$ is nilpotent, then the following are equivalent:
\begin{enumerate}[label=\textup{(\roman*)}]
\item $\cC(T,\Gamma)=0$;
\item $K=0$;
\item the unital regular representation $\widehat\rho_T$ is faithful.
\end{enumerate}
Thus a nonzero nilpotent operator defect is always detected in first relative homology.
\end{proposition}

\begin{proof}
Only (i)$\Rightarrow$(ii) requires proof.  The equality $K/K^2=0$ gives $K=K^2$.  Iteration yields $K=K^{2^m}$ for every $m$.  If $K^N=0$, choose $2^m\ge N$ to obtain $K=0$.
\end{proof}

\begin{remark}
Nilpotence is essential.  For a binary-derived unital ring, the Dorroh kernel is idempotent: $K=K^2\ne0$, so the conormal defect vanishes although the unitalization map is not injective.  The theorem therefore separates finite-depth invisible structure from identity-normalization artifacts of infinite depth.
\end{remark}

\begin{proposition}[Two extreme towers]\label{prop:extreme-towers}
For every two-step tensor system of \cref{sec:tensor-systems}, $K^2=0$, hence $\cC_n=0$ and $\cS_n=\cC^{\otimes n}$ for $n\ge2$.  For the rank-one family $T_d$ of \cref{sec:rank-one}, every map $\tau_n$ is an isomorphism; hence its infinite defect tower is tensor-free.
\end{proposition}

\begin{proof}
The first assertion follows from \cref{thm:tensor-structure}.  For $T_d$, the class of $\kappa^{n}=(-d^2)^{n-1}\kappa$ generates $K^n/K^{n+1}$.  Since $\cC\cong\ZZ/d^2\ZZ$, its $n$-fold tensor power over $\ZZ$ is again $\ZZ/d^2\ZZ$, and multiplication sends a generator to a generator up to sign.
\end{proof}

\section{Square-zero regular reductions and Hochschild curvature}

Assume in this section that $K^2=0$ and that the quotient in \eqref{eq:quotient} admits a unital additive section, meaning an additive section that carries $1_{\cR}$ to $1_{\cA}$.  This is slightly stronger than a bare splitting of the underlying abelian-group sequence.  It holds for all finite-free examples below: there $\cR\cong\ZZ\oplus I$ additively, with the identity in the displayed $\ZZ$-summand.

Here $\HH^2_{\ZZ}(\cR,K)$ is computed from the normalized
$\ZZ$-multilinear Hochschild cochains
$C^n=\Hom_{\ZZ}(\cR^{\otimes_{\ZZ}n},K)$ with the usual Hochschild
differential; normalized cochains vanish whenever one input is $1_{\cR}$.
Two square-zero extensions are equivalent when there is a unital ring
isomorphism between them inducing the identity on both the ideal $K$ and the
quotient $\cR$.  The extension interpretation used below classifies precisely
the unital square-zero extensions that split as $\ZZ$-modules; it does not
classify extensions whose underlying abelian-group sequence is nonsplit.

Choose a unital additive section $s:\cR\to\cA$. Define
\begin{equation}\label{eq:cocycle}
\omega_s(r,r')=s(r)s(r')-s(rr')\in K.
\end{equation}
The ideal $K$ is naturally an $\cR$-bimodule because $K^2=0$.

\begin{proposition}[Regular extension class]\label{thm:HHclass}
The map $\omega_s$ is a normalized Hochschild $2$-cocycle. Its cohomology class
\[
\Omega(T,\Gamma)=[\omega_s]\in\HH^2_{\ZZ}(\cR_T,K)
\]
is independent of the additive section $s$. Moreover, the following are equivalent:
\begin{enumerate}[label=\textup{(\roman*)}]
\item $\Omega(T,\Gamma)=0$;
\item the quotient $\cA(T,\Gamma)\twoheadrightarrow\cR_T$ admits a unital ring section;
\item the extension is equivalent to the split square-zero extension $\cR_T\ltimes K$.
\end{enumerate}
\end{proposition}

\begin{proof}
Associativity of multiplication in $\cA$ gives
\[
r\omega_s(r',r'')-\omega_s(rr',r'')
+\omega_s(r,r'r'')-\omega_s(r,r')r''=0,
\]
which is the Hochschild cocycle identity. If $s'=s+h$ for an additive normalized map $h:\cR\to K$, then $K^2=0$ gives
\[
\omega_{s'}-\omega_s=\delta h.
\]
Thus the class is independent of $s$. If the class vanishes, choose $h$ with $\omega_s=\delta h$; then $s-h$ is multiplicative. The converse is immediate. The equivalence with a split square-zero extension is the standard extension interpretation of $\HH^2$ \cite{Hochschild1945,Witherspoon2019}.
\end{proof}

\begin{criterion}[A nonsplitting test]\label{crit:nonsplit}
Suppose $K^2=0$. Let $v\in\cR$ satisfy $v^2=0$, and let $p\in\cA$ be a lift of $v$. If
\[
p^2\ne0\quad\text{and}\quad pK=Kp=0,
\]
then $\Omega(T,\Gamma)\ne0$.
\end{criterion}

\begin{proof}
Every lift of $v$ has the form $p+k$. Its square is $p^2+pk+kp+k^2=p^2\ne0$. A ring section would have to send $v$ to a square-zero lift.
\end{proof}

\section{Two-step nilpotent tensor systems}\label{sec:tensor-systems}

We now describe a large class for which the operator envelope, conormal defect, extension class and characteristic behavior reduce to exact linear algebra.

Let $E$ and $W$ be nonzero finite-rank free abelian groups and let
\[
\mu:E\otimes E\otimes E\longrightarrow W
\]
be a surjective trilinear map.

\begin{construction}\label{con:Tmu}
Set $T_\mu=E\oplus W$, take $\Gamma=\ZZ$, and define
\[
[x+u,y+v,z+w]=\mu(x,y,z),
\qquad x,y,z\in E,\quad u,v,w\in W.
\]
Equivalently,
\[
(x+u)\,n\,(y+v)\,m\,(z+w)=nm\mu(x,y,z).
\]
\end{construction}

\begin{proposition}\label{prop:Tmu-assoc}
The construction $T_\mu$ is an additive ternary $\ZZ$-ring. Every iterated ternary product is zero.
\end{proposition}

\begin{proof}
Every product lies in $W$, and the defining product vanishes whenever any input lies in $W$. Hence all three iterated products in ternary associativity are zero.
\end{proof}

Put
\[
P=E\otimes E,
\quad L=W\otimes E,
\quad R=E\otimes W.
\]
Define
\begin{align}
\Theta_\mu:P&\longrightarrow\Hom_{\ZZ}(E,W),
&\Theta_\mu(a\otimes b)(c)&=\mu(a,b,c),\label{eq:Theta}\\
B_\mu:E^{\otimes4}&\longrightarrow L\oplus R,
&B_\mu(a\otimes b\otimes c\otimes d)
&=\mu(a,b,c)\otimes d-a\otimes\mu(b,c,d).
\label{eq:Bmu}
\end{align}
Let
\[
Q_\mu=\coker B_\mu.
\]
Finally, define a bilinear map
\begin{equation}\label{eq:mbar}
\overline m_\mu:P\times P\to Q_\mu,
\qquad
\overline m_\mu(a\otimes b,c\otimes d)
=\overline{\mu(a,b,c)\otimes d}.
\end{equation}

\begin{theorem}[Operator structure theorem]\label{thm:tensor-structure}
For the tensor system $T_\mu$ there is an isomorphism of abelian groups
\[
\cO(T_\mu,\ZZ)\cong P\oplus Q_\mu
\]
under which multiplication is
\begin{equation}\label{eq:tensor-mult}
(p,q)(p',q')=(0,\overline m_\mu(p,p')).
\end{equation}
Consequently,
\[
\cO(T_\mu,\ZZ)^2\subseteq Q_\mu,
\qquad
\cO(T_\mu,\ZZ)^3=0.
\]
The regular representation is $\Theta_\mu$ on $P$ and zero on $Q_\mu$. Hence
\begin{equation}\label{eq:D-tensor}
\Def(T_\mu,\ZZ)
\cong\Ker\Theta_\mu\oplus Q_\mu.
\end{equation}
Furthermore,
\[
\UDef(T_\mu,\ZZ)=\Def(T_\mu,\ZZ),
\qquad
\UDef(T_\mu,\ZZ)^2=0,
\]
and
\[
\cR_{T_\mu}\cong\ZZ\ltimes\im\Theta_\mu
\]
with square-zero ideal $\im\Theta_\mu$.
\end{theorem}

\begin{proof}
Since $\Gamma=\ZZ$, the tensor group underlying the raw operator rng is $T_\mu\otimes T_\mu$, which decomposes as
\[
P\oplus L\oplus R\oplus(W\otimes W).
\]
The middle relations with four entries in $E$ are exactly the image of $B_\mu$. Relations with three entries in $E$ and the remaining outer entry in $W$ kill $W\otimes W$ because $\mu$ is surjective. No middle relation has a component in $P$. This gives the stated additive decomposition.

Products of two elements of $P$ are given by \eqref{eq:mbar}. A product involving a representative in $L$ or $R$ vanishes: either one of the first three ternary inputs lies in $W$, or the result lies in the already killed summand $W\otimes W$. Hence \eqref{eq:tensor-mult} and cubic nilpotence follow.

On the regular module, $a\otimes b$ acts as $c+u\mapsto\mu(a,b,c)$, while $L$ and $R$ act trivially. This proves \eqref{eq:D-tensor}. If an operator in the image were a nonzero integer multiple of $\id_{T_\mu}$, restricting to the nonzero torsion-free summand $W$ would give a contradiction. Thus $I_T=0$ in \cref{prop:def-seq}, so the intrinsic and unital defects coincide.

If $k\in\Ker\Theta_\mu$, then for every $c\otimes d\in P$,
\[
k(c\otimes d)=\Theta_\mu(k)(c)\otimes d=0.
\]
Elements of $Q_\mu$ annihilate the whole operator rng on the left. Therefore the square of the defect is zero. Finally, compositions of visible regular operators vanish because they map $E$ into $W$ and kill $W$.
\end{proof}

\begin{remark}[Scalar form of the structure theorem]
\label{rem:tensor-structure-base}
The presentation argument is base-linear.  For every commutative unital
$\ZZ$-algebra $S$, the split surjection $\mu$ remains surjective after scalar
extension, and the same proof gives
\[
\cO_S(T_{\mu,S})\cong P_S\oplus Q_{\mu,S},
\qquad
\Def_S(T_{\mu,S})
\cong\Ker\Theta_{\mu,S}\oplus Q_{\mu,S},
\]
as $S$-modules, where $P_S=P\otimes_{\ZZ}S$ and
$Q_{\mu,S}=\coker B_{\mu,S}$; multiplication is given by the scalar-extended
version of \eqref{eq:tensor-mult}.  In particular, the dimension formulas used
below over fields are instances of this scalar form, not an implicit
identification with the integral kernel.
\end{remark}

\begin{corollary}[Rank and torsion profile]\label{cor:rank-profile}
Let $r=\rank E$ and $s=\rank W$. Then
\begin{align*}
\rank\Def(T_\mu,\ZZ)
&=r^2-\rank\Theta_\mu+2rs-\rank B_\mu,\\
\operatorname{tors}\Def(T_\mu,\ZZ)
&\cong\operatorname{tors}(\coker B_\mu).
\end{align*}
The same statements hold for the conormal defect.
\end{corollary}

\begin{proof}
Use \eqref{eq:D-tensor}, the fact that $\Ker\Theta_\mu$ is free, and $\UDef^2=0$.
\end{proof}

\begin{definition}\label{def:curvature}
Choose an additive splitting $\sigma:\im\Theta_\mu\to P$. The \emph{operator curvature cocycle} is
\[
\omega_\mu(v,w)=\overline m_\mu(\sigma(v),\sigma(w))
\in Q_\mu\subseteq\Def(T_\mu,\ZZ).
\]
Its Hochschild class is the regular extension class $\Omega(T_\mu,\ZZ)$.
\end{definition}

This gives a concrete tensor formula for the abstract class of \cref{thm:HHclass}.

\section{Smith normal form and characteristic jumps}\label{sec:snf}

Fix bases $e_1,\dots,e_r$ of $E$ and $w_1,\dots,w_s$ of $W$, and write
\[
\mu(e_a,e_b,e_c)=\sum_{u=1}^s\mu^u_{abc}w_u.
\]
Then $\Theta_\mu$ is the $(sr)\times r^2$ integer matrix
\begin{equation}\label{eq:Theta-matrix}
(\Theta_\mu)_{(u,c),(a,b)}=\mu^u_{abc},
\end{equation}
and $B_\mu$ is the $(2sr)\times r^4$ integer matrix whose column $(a,b,c,d)$ has entries
\begin{equation}\label{eq:B-matrix}
\mu^u_{abc}\quad\text{in the row }(L;u,d),
\qquad
-\mu^u_{bcd}\quad\text{in the row }(R;a,u).
\end{equation}

\begin{proposition}[Exact arithmetic extraction]\label{thm:SNF}
The Smith normal form of $B_\mu$ determines the complete torsion subgroup and free rank of $Q_\mu$. A $\ZZ$-basis of $\Ker\Theta_\mu$, together with Smith generators for $Q_\mu$, determines the conormal defect. The multiplication table and the curvature cocycle are obtained by reducing the vectors
\[
\mu(a,b,c)\otimes d
\]
modulo the Smith relations of $B_\mu$.
\end{proposition}

\begin{proof}
This is the standard classification of finitely generated abelian groups applied to $Q_\mu=\coker B_\mu$, together with the explicit structure theorem \cref{thm:tensor-structure}.
\end{proof}

\begin{definition}\label{def:arith-pair}
The \emph{visible saturation defect} and the \emph{middle-relation torsion defect} are
\[
\mathsf S_{\mathrm{vis}}(\mu)
=\operatorname{tors}(\coker\Theta_\mu),
\qquad
\mathsf S_{\mathrm{mid}}(\mu)
=\operatorname{tors}(\coker B_\mu).
\]
The ordered pair
\[
\mathsf S_{\mathrm{arith}}(\mu)
=(\mathsf S_{\mathrm{vis}}(\mu),\mathsf S_{\mathrm{mid}}(\mu))
\]
is called the \emph{arithmetic operator-defect pair}.
\end{definition}

The first group measures the failure of the visible operator lattice to be saturated inside $\Hom(E,W)$. The second is the torsion already present in the hidden middle-relation quotient $Q_\mu$. They are logically independent.

\begin{proposition}[Smith decomposition of the jump]\label{prop:smith-jump}
Let
\[
\alpha_1\mid\cdots\mid\alpha_{a_0}
\quad\text{and}\quad
\beta_1\mid\cdots\mid\beta_{b_0}
\]
be the nonzero Smith invariant factors of $\Theta_\mu$ and $B_\mu$, respectively. For every prime $p$ for which $\mu_p$ remains surjective,
\begin{align*}
a_0-a_p&=\#\{i:p\mid\alpha_i\},\\
b_0-b_p&=\#\{j:p\mid\beta_j\}.
\end{align*}
Consequently, the characteristic jump equals the number of $p$-divisible invariant factors occurring in the two components of $\mathsf S_{\mathrm{arith}}(\mu)$.
\end{proposition}

\begin{proof}
After Smith reduction, an invariant factor remains nonzero modulo $p$ exactly
when it is not divisible by $p$.  Thus the two rank losses are the displayed
cardinalities.  Over a field $k$, \eqref{eq:D-tensor} gives
\[
\dim_k\Def_k(T_{\mu_k})
=r^2-\rank_k\Theta_{\mu,k}+2rs-\rank_k B_{\mu,k}.
\]
Subtracting the formula over $\QQ$ from the formula over $\FF_p$ shows directly
that the total jump is the sum of those two rank losses.
\end{proof}

\begin{proposition}[Isomorphism invariance]\label{prop:isom-invariance}
If $g:E\to E'$ and $h:W\to W'$ are unimodular isomorphisms carrying $\mu$ to $\mu'$, then the matrices $\Theta_\mu,B_\mu$ and $\Theta_{\mu'},B_{\mu'}$ differ by unimodular changes of rows and columns. Hence the free defect rank, the arithmetic operator-defect pair, the bad-prime set and the regular extension class are invariants of the integral ternary system.
\end{proposition}

\begin{proof}
The maps $g^{\otimes2}$, $g^{\otimes4}$, $h\otimes g$, $g\otimes h$ and the induced conjugation on $\Hom(E,W)$ intertwine the definitions of $\Theta$ and $B$. In chosen bases these maps are unimodular. The extension class is functorial under the induced isomorphisms of the operator and regular rings.

It remains to justify that this change-of-bases formulation captures an
arbitrary isomorphism of the ternary systems.  Surjectivity of $\mu$ implies
that $W$ is exactly the subgroup of $T_\mu$ generated by all ternary products.
It is therefore characteristic under ternary-ring isomorphisms.  An
isomorphism $T_\mu\to T_{\mu'}$ restricts to a unimodular map
$h:W\to W'$ and induces a unimodular map
$g:T_\mu/W\to T_{\mu'}/W'$, identified with $E\to E'$.  Relative to chosen
complements, the isomorphism may also have a shear $E\to W'$, but every
ternary product ignores $W'$-components in its inputs.  Hence the shear drops
out and $h\mu=\mu' g^{\otimes3}$.  The preceding intertwining argument now
applies, proving invariance for the integral ternary system itself.
\end{proof}

\section{Independent realization of arithmetic defects}\label{sec:realization}

The matrices $\Theta_\mu$ and $B_\mu$ are built from the same tensor, so their torsion cokernels might appear constrained by one another.  We now prove that no such constraint exists at the level of finite abelian groups.  The proof is constructive and uses two complementary atomic families.

\begin{construction}[Visible and middle atoms]\label{con:atoms}
For $n\ge1$, let $E=\ZZ e_0\oplus\ZZ e_1$ and $W=\ZZ w$.  Define surjective tensors
\(V_n,M_n:E^{\otimes3}\to W\) by the following table, whose columns are ordered lexicographically:
\[
\begin{array}{c|rrrrrrrr}
(a,b,c)&000&001&010&011&100&101&110&111\\ \hline
V_n(e_a,e_b,e_c)/w&-n&0&n&n&-1&n&-1&0\\
M_n(e_a,e_b,e_c)/w&0&n&0&0&0&-n&0&1
\end{array}
\]
Surjectivity follows from the displayed coefficients $-1$ and $1$, respectively.
\end{construction}

\begin{theorem}[Atomic separation]\label{thm:atomic-separation}
For every $n\ge1$ the two atoms have complementary arithmetic defect pairs:
\[
\mathsf S_{\mathrm{arith}}(V_n)
   =(\ZZ/n\ZZ,0),
\qquad
\mathsf S_{\mathrm{arith}}(M_n)
   =(0,\ZZ/n\ZZ).
\]
More precisely,
\begin{align*}
\operatorname{SNF}(\Theta_{V_n})&=\operatorname{diag}(1,n),
&\operatorname{SNF}(B_{V_n})&=\operatorname{diag}(1,1,1,1),\\
\operatorname{SNF}(\Theta_{M_n})&=\operatorname{diag}(1,0),
&\operatorname{SNF}(B_{M_n})&=\operatorname{diag}(1,1,1,n).
\end{align*}
Here the zero in the Smith form of $\Theta_{M_n}$ records a free cokernel summand and contributes no torsion.
\end{theorem}

\begin{proof}
With columns indexed by $(00,01,10,11)$, direct substitution gives
\[
\Theta_{V_n}=
\begin{pmatrix}
-n&n&-1&-1\\
0&n&n&0
\end{pmatrix},
\qquad
\Theta_{M_n}=
\begin{pmatrix}
0&0&0&0\\
n&0&-n&1
\end{pmatrix}.
\]
For $\Theta_{V_n}$ the gcd of the entries is $1$, every $2\times2$ minor is divisible by $n$, and the minor on the last two columns is $n$.  Its determinantal divisors are therefore $1,n$.  The nonzero row of $\Theta_{M_n}$ has relatively prime entries, giving the stated Smith form.

Index the columns of each middle-relation matrix by the binary words $abcd$ and
the rows by $(L_0,L_1,R_0,R_1)$.  To make the all-$n$ calculation explicit,
split the columns into the eight words beginning with $0$ and the eight
beginning with $1$, in lexicographic order.  Direct substitution in
\eqref{eq:B-matrix} gives
{\setlength{\arraycolsep}{3pt}
\[
\begin{aligned}
B_{V_n}^{[0]}&=
\begin{pmatrix}
-n&0&0&0&n&0&n&0\\
0&-n&0&0&0&n&0&n\\
n&0&-n&-n&1&-n&1&0\\
0&0&0&0&0&0&0&0
\end{pmatrix},&
B_{V_n}^{[1]}&=
\begin{pmatrix}
-1&0&n&0&-1&0&0&0\\
0&-1&0&n&0&-1&0&0\\
0&0&0&0&0&0&0&0\\
n&0&-n&-n&1&-n&1&0
\end{pmatrix},\\[1ex]
B_{M_n}^{[0]}&=
\begin{pmatrix}
0&0&n&0&0&0&0&0\\
0&0&0&n&0&0&0&0\\
0&-n&0&0&0&n&0&-1\\
0&0&0&0&0&0&0&0
\end{pmatrix},&
B_{M_n}^{[1]}&=
\begin{pmatrix}
0&0&-n&0&0&0&1&0\\
0&0&0&-n&0&0&0&1\\
0&0&0&0&0&0&0&0\\
0&-n&0&0&0&n&0&-1
\end{pmatrix}.
\end{aligned}
\]
}
Here $B^{[0]}$ has columns $0000,\ldots,0111$ and $B^{[1]}$ has columns
$1000,\ldots,1111$.

For $B_{V_n}$, the full minor on columns
$0100,1000,1001,1110$ is
\[
\det\begin{pmatrix}
n&-1&0&0\\0&0&-1&0\\1&0&0&0\\0&n&0&1
\end{pmatrix}=1,
\]
so $B_{V_n}$ is surjective over $\ZZ$.

For $B_{M_n}$ write $C_w$ for the column indexed by $w$ and put
\[
u_1=C_{0111}=(0,0,-1,0)^{\mathsf T},\quad
u_2=C_{1110}=(1,0,0,0)^{\mathsf T},\quad
u_3=C_{1111}=(0,1,0,-1)^{\mathsf T}.
\]
The displayed matrices give the complete reduction
\[
C_w\equiv0\pmod n
\quad(w\notin\{0111,1110,1111\}),
\qquad
(C_{0111},C_{1110},C_{1111})\equiv(u_1,u_2,u_3)\pmod n.
\]
Consequently every four columns are dependent modulo $n$, so every full minor
is divisible by $n$.  On rows $(L_0,L_1,R_0)$ and columns
$0111,1110,1111$ the determinant is $-1$; hence the first three
determinantal divisors are $1$.  Finally, on columns
$0011,0111,1110,1111$ the full determinant is
\[
\det\begin{pmatrix}
0&0&1&0\\ n&0&0&1\\0&-1&0&0\\0&0&0&-1
\end{pmatrix}=n.
\]
Thus the fourth determinantal divisor is exactly $n$.  This proves both Smith
profiles, including the quantified statement for every $n\ge1$, and hence the
asserted torsion cokernels.
\end{proof}

\begin{proposition}[Orthogonal-sum theorem]\label{prop:orthogonal-sum}
Let $\mu_i:E_i^{\otimes3}\twoheadrightarrow W_i$ be surjective tensors of finite-rank free abelian groups, for $1\le i\le t$.  Put
\[
E=\bigoplus_{i=1}^tE_i,
\qquad
W=\bigoplus_{i=1}^tW_i,
\qquad
\mu(x,y,z)=\bigoplus_{i=1}^t
\mu_i(x_i,y_i,z_i).
\]
Then $\mu:E^{\otimes3}\twoheadrightarrow W$ is surjective and
\begin{align}
\coker\Theta_\mu
&\cong
\left(\bigoplus_{i=1}^t\coker\Theta_{\mu_i}\right)
\oplus
\left(\bigoplus_{i\ne j}\Hom_{\ZZ}(E_j,W_i)\right),
\label{eq:theta-orthogonal}\\
\coker B_\mu
&\cong\bigoplus_{i=1}^t\coker B_{\mu_i}.
\label{eq:B-orthogonal}
\end{align}
Consequently,
\[
\mathsf S_{\mathrm{arith}}(\mu)
=
\left(
\bigoplus_i\mathsf S_{\mathrm{vis}}(\mu_i),
\bigoplus_i\mathsf S_{\mathrm{mid}}(\mu_i)
\right).
\]
\end{proposition}

\begin{proof}
The flattening $\Theta_\mu$ is the direct sum of the diagonal maps $\Theta_{\mu_i}$.  It is zero on $E_i\otimes E_j$ for $i\ne j$, and the off-diagonal target blocks $\Hom(E_j,W_i)$ are free.  This proves \eqref{eq:theta-orthogonal}.

For $B_\mu$, the summand with all four entries in $E_i$ is exactly $B_{\mu_i}$.  If the first three entries lie in $E_i$ and the fourth lies in $E_j$ with $i\ne j$, the corresponding relations surject onto $W_i\otimes E_j$, because $\mu_i$ is surjective.  Similarly, quadruples with the last three entries in $E_j$ kill $E_i\otimes W_j$.  These are all off-diagonal target summands, and no relation couples them to a diagonal summand.  Hence they contribute zero cokernel and \eqref{eq:B-orthogonal} follows.
\end{proof}

\begin{theorem}[Independent arithmetic realization]\label{thm:independent-realization}
For every ordered pair $(G,H)$ of finite abelian groups, there exist finite-rank free abelian groups $E,W$ and a surjective tensor $\mu:E^{\otimes3}\to W$ such that the torsion-free two-step ternary ring $T_\mu=E\oplus W$ satisfies
\[
\mathsf S_{\mathrm{arith}}(\mu)=(G,H).
\]
If
\[
G\cong\bigoplus_{i=1}^{a}\ZZ/n_i\ZZ,
\qquad
H\cong\bigoplus_{j=1}^{b}\ZZ/m_j\ZZ
\quad(n_i,m_j\ge2),
\]
and $a+b>0$, one may take
\[
\rank E=2(a+b),
\qquad
\rank W=a+b.
\]
Thus the invariant factors and prime supports of the visible and middle-relation defects can be prescribed independently.
\end{theorem}

\begin{proof}
Take the orthogonal sum of the visible atoms $V_{n_i}$ and the middle atoms $M_{m_j}$.  The tensor is surjective and the underlying group $E\oplus W$ is free abelian.  The atomic Smith profiles in \cref{thm:atomic-separation} and the decomposition in \cref{prop:orthogonal-sum} give exactly $(G,H)$.  The displayed ranks follow because every atom has $\rank E=2$ and $\rank W=1$.  If $G=H=0$, the atom $V_1$ supplies a nonzero torsion-free realization with trivial arithmetic pair.
\end{proof}

\begin{corollary}[Independent bad-prime supports]\label{cor:independent-primes}
For arbitrary finite sets of primes $P$ and $Q$, there is a torsion-free two-step ternary $\ZZ$-ring whose visible saturation defect is supported exactly on $P$ and whose middle-relation torsion defect is supported exactly on $Q$.
\end{corollary}

\begin{proof}
Apply \cref{thm:independent-realization} to
$G=\bigoplus_{p\in P}\ZZ/p\ZZ$ and
$H=\bigoplus_{q\in Q}\ZZ/q\ZZ$.
\end{proof}

\begin{remark}
The classification theorem for finite abelian groups is used only to assemble the cyclic atoms.  The complementary Smith profiles and the orthogonal decomposition are properties of the ternary presentation matrices; they do not follow from the existence of an associative envelope or from ordinary module theory.
\end{remark}

Let $k$ be a field and let $\mu_k$ denote reduction or scalar extension. Assume $\mu_k$ remains surjective.

\begin{proposition}[Characteristic defect jump]\label{prop:char-jump}
Let
\[
a_0=\rank_{\QQ}\Theta_\mu,
\qquad b_0=\rank_{\QQ}B_\mu,
\]
and, for a prime $p$,
\[
a_p=\rank_{\FF_p}(\Theta_\mu\bmod p),
\qquad b_p=\rank_{\FF_p}(B_\mu\bmod p).
\]
Then
\begin{align}
\dim_{\FF_p}\Def_{\FF_p}(T_{\mu_p})
&=r^2-a_p+2rs-b_p,\label{eq:dim-p}\\
\dim_{\FF_p}\Def_{\FF_p}(T_{\mu_p})
-\dim_{\QQ}\Def_{\QQ}(T_{\mu_\QQ})
&=(a_0-a_p)+(b_0-b_p).\label{eq:jump}
\end{align}
Only finitely many primes have a nonzero jump.
\end{proposition}

\begin{proof}
Apply the scalar form of \cref{thm:tensor-structure} recorded in
\cref{rem:tensor-structure-base}. Matrix ranks can drop modulo only those
primes dividing suitable nonzero minors of $\Theta_\mu$ or $B_\mu$.
\end{proof}

\begin{remark}\label{rem:nonflat}
The actual defect after reduction modulo $p$ need not equal the integral defect tensored with $\FF_p$. Torsion in $Q_\mu$ contributes one source of extra dimension, while a rank drop of $\Theta_\mu$ contributes a second source. The sign-tensor example in \cref{sec:sign} exhibits both simultaneously.
\end{remark}

\section{Base change, operator discriminants and defect strata}\label{sec:basechange}

By \cref{thm:envelope-basechange}, the universal operator envelope commutes with scalar extension.  The integral tensor theorem now separates the part of the \emph{regular kernel} that commutes with scalar extension from the part created by nonsaturated visible operators.
The Tor correction, determinantal rank loci and Fitting-ideal base-change law
used in this section are standard.  The specific result is the two-matrix
decomposition: it identifies the correction term belonging to regular
visibility while showing that the middle-relation cokernel is already inherited
by scalar extension.

Let $C_{\mathrm{vis}}(\mu)=\coker\Theta_\mu$.  For a commutative unital
$\ZZ$-algebra $S$, write a subscript $S$ for tensoring the entire two-step
system with $S$.

\begin{theorem}[Base-change obstruction exact sequence]\label{thm:basechange}
For every commutative unital $\ZZ$-algebra $S$ there is a canonical exact sequence
\begin{equation}\label{eq:basechange-exact}
0\longrightarrow
\Def(T_\mu,\ZZ)\otimes_{\ZZ}S
\longrightarrow
\Def_S(T_{\mu,S})
\longrightarrow
\Tor^{\ZZ}_1(C_{\mathrm{vis}}(\mu),S)
\longrightarrow0.
\end{equation}
Consequently the operator defect commutes with every flat base change.  The middle-relation quotient always satisfies
\(Q_\mu\otimes S\cong Q_{\mu,S}\); all failure of naive defect base change is controlled by the visible saturation cokernel.
\end{theorem}

\begin{proof}
Put $P=E\otimes_{\ZZ}E$ and $H=\Hom_{\ZZ}(E,W)$, and factor
$\Theta_\mu:P\to H$ through its image $I$.  Because $E$ and $W$ are finite
free, there are canonical identifications
\[
P\otimes_{\ZZ}S\cong E_S\otimes_S E_S,
\qquad
H\otimes_{\ZZ}S\cong\Hom_S(E_S,W_S),
\]
under which $\Theta_\mu\otimes1$ is exactly the scalar-extended flattening
$\Theta_{\mu,S}$.  The groups $P,H,I$ are free abelian.  In particular, the
sequence
$0\to\Ker\Theta_\mu\to P\to I\to0$ splits because $I$ is projective, so it
remains a short exact sequence after tensoring with arbitrary $S$.
Tensoring $0\to I\to H\to C_{\mathrm{vis}}\to0$ gives
\[
0\to\Tor_1^{\ZZ}(C_{\mathrm{vis}},S)
\to I\otimes S\to H\otimes S.
\]
Taking the inverse image in $P\otimes S$ yields
\[
0\to(\Ker\Theta_\mu)\otimes S
\to\Ker(\Theta_\mu\otimes S)
\to\Tor_1^{\ZZ}(C_{\mathrm{vis}},S)\to0.
\]
Since the surjection $\mu:E^{\otimes3}\to W$ splits as a map of abelian
groups, its scalar extension remains surjective.  The scalar form in
\cref{rem:tensor-structure-base} therefore gives
\[
\Def_S(T_{\mu,S})
\cong\Ker\Theta_{\mu,S}\oplus Q_{\mu,S}.
\]
Because cokernels commute with tensor products,
$Q_\mu\otimes S\cong\coker(B_\mu\otimes S)=Q_{\mu,S}$.  Taking the identity
on this common summand and the preceding exact sequence on the visible summand
gives \eqref{eq:basechange-exact}.
\end{proof}

\begin{corollary}[Prime-fiber correction term]\label{cor:prime-fiber}
For every prime $p$,
\[
0\to\Def(T_\mu,\ZZ)\otimes\FF_p
\to\Def_{\FF_p}(T_{\mu_p})
\to C_{\mathrm{vis}}(\mu)[p]\to0.
\]
Hence the excess dimension over the tensor-reduced integral defect is exactly the number of Smith invariant factors of $\Theta_\mu$ divisible by $p$.
\end{corollary}

\begin{definition}[Operator discriminant]\label{def:opdisc}
Let $a_0=\rank_{\QQ}\Theta_\mu$ and $b_0=\rank_{\QQ}B_\mu$.  If $a_0>0$, let $\Delta_{\mathrm{vis}}(\mu)$ be the positive gcd of all nonzero $a_0\times a_0$ minors of $\Theta_\mu$, and put it equal to $1$ if $a_0=0$.  Define $\Delta_{\mathrm{mid}}(\mu)$ analogously from $B_\mu$, and set
\[
\Delta_{\mathrm{op}}(\mu)
=\Delta_{\mathrm{vis}}(\mu)\Delta_{\mathrm{mid}}(\mu).
\]
\end{definition}

\begin{theorem}[Arithmetic discriminant criterion]\label{thm:disc-criterion}
A prime $p$ has a nonzero characteristic defect jump if and only if
\(p\mid\Delta_{\mathrm{op}}(\mu)\).  More precisely, the jump is
\[
\#\{i:p\mid\alpha_i\}
+\#\{j:p\mid\beta_j\},
\]
where the $\alpha_i$ and $\beta_j$ are the nonzero Smith invariant factors of $\Theta_\mu$ and $B_\mu$.  Thus $\Delta_{\mathrm{op}}$ cuts out the complete bad-prime support, while the Smith factors retain its multiplicities.
\end{theorem}

\begin{proof}
A matrix loses rank modulo $p$ exactly when all of its maximal-rank minors vanish modulo $p$.  Over $\ZZ$, this is equivalent to divisibility of their gcd by $p$.  The quantitative statement is \cref{prop:smith-jump}.
\end{proof}

\subsection{Geometric parameter strata}
Let $k$ be a field, let $E$ and $W$ have dimensions $r$ and $s$, and let
\[
\cV=\Hom_k(E^{\otimes3},W).
\]
The surjective tensors form a Zariski-open subset $\cV^{\mathrm{surj}}$.  The matrices $\Theta_\mu$ and $B_\mu$ have entries linear in the coordinates of $\mu$.

\begin{proposition}[Determinantal defect stratification]\label{thm:defect-strata}
For integers $a,b$, the subsets
\[
\cX_{a,b}
=\{\mu\in\cV^{\mathrm{surj}}:
\rank\Theta_\mu\le a,
\ \rank B_\mu\le b\}
\]
are closed in $\cV^{\mathrm{surj}}$ and stable under the natural action of $\mathrm{GL}(E)\times\mathrm{GL}(W)$.  The function
\[
\delta(\mu)=\dim_k\Def_k(T_\mu)
=r^2-\rank\Theta_\mu+2rs-\rank B_\mu
\]
is upper semicontinuous.  Consequently, every irreducible component of $\cV^{\mathrm{surj}}$ contains a dense open subset on which the defect dimension is minimal and constant.
\end{proposition}

\begin{proof}
The conditions on the two ranks are the vanishing of minors, hence define determinantal closed subsets.  Equivariance of $\Theta$ and $B$ under change of bases proves stability.  Matrix rank is lower semicontinuous, so each corank and their sum $\delta$ are upper semicontinuous.  An upper semicontinuous integer-valued function assumes its minimum on a dense open subset of each irreducible component.
\end{proof}

\begin{remark}[Fitting formulation]
The same strata are encoded scheme-theoretically by the Fitting ideals of $\coker\Theta$ and $\coker B$.  Since Fitting ideals commute with base change \cite{StacksFitting}, they provide a coordinate-free version of the arithmetic and geometric defect loci.  In particular, $\Delta_{\mathrm{op}}$ records the bad-prime support of these maximal-rank Fitting strata; no new general discriminant theory is asserted.
\end{remark}

\subsection{Directional visibility and triality}
The regular left-module envelope privileges the flattening $(12|3)$.  For a tensor $\mu$, let $\Theta^{12|3}_\mu$, $\Theta^{23|1}_\mu$ and $\Theta^{13|2}_\mu$ denote the three flattenings obtained by selecting two input slots as operator slots and the remaining slot as the acted-on variable.

\begin{definition}
For finite-free abelian groups $E,W$, the \emph{integral triality visibility
profile} of $\mu$ is the unordered triple of Smith profiles of the three
flattenings.  For finite-dimensional vector spaces over a field, the
\emph{field triality visibility profile} is the unordered triple of their
ranks.
\end{definition}

\begin{proposition}\label{prop:triality}
Over $\ZZ$, the integral triality visibility profile is invariant under
$\mathrm{GL}(E)\times\mathrm{GL}(W)$; over a field, the field triality
visibility profile has the analogous invariance.  In either setting the
profile is permuted by the action of the symmetric group on the three input
slots.  If $\mu$ is symmetric or alternating, all three constituent profiles
coincide.
\end{proposition}

\begin{proof}
A change of bases multiplies each flattening by invertible matrices on the source and target.  Permuting tensor slots permutes the flattenings.  Symmetry or alternation identifies them up to signs and permutation matrices, which do not alter Smith factors or ranks.
\end{proof}

\section{A non-binary example with all-degree periodic homology}\label{sec:nil-example}

The periodic mechanism below is classical: it is the rank-one matrix
factorization $(p,p^2)$ of the hypersurface equation $p^3$, in the sense of
Eisenbud \cite{Eisenbud1980}.  The contribution here is its realization as the
universal-to-regular operator quotient of a ternary ring that is not induced by
an associative binary product.  In this section, \emph{non-binary} means
precisely ``not induced by an associative $\ZZ$-bilinear binary multiplication
on the underlying abelian group.''

Let $T=\ZZ e\oplus\ZZ f$, let $\Gamma=\ZZ$, and define
\[
[e,e,e]=f,
\]
with every other basis product zero.

\begin{proposition}\label{prop:nil-assoc}
This defines an additive ternary $\ZZ$-ring, and it is not induced by any associative $\ZZ$-bilinear binary multiplication on $T$.
\end{proposition}

\begin{proof}
All iterated ternary products vanish, so associativity is immediate. Suppose $[x,y,z]=(x\star y)\star z$ for an associative bilinear product. After tensoring with $\QQ$, write
\[
e\star e=ae+bf,
\qquad f\star e=ge+hf.
\]
The identities $(e\star e)\star e=f$ and $(f\star e)\star e=0$ imply
\[
a^2+bg=0,
\quad b(a+h)=1,
\quad g(a+h)=0,
\quad gb+h^2=0.
\]
Since $a+h\ne0$, one gets $g=0$, then $a=h=0$, contradicting $b(a+h)=1$.
\end{proof}

\begin{theorem}\label{thm:nil-envelope}
For this ternary ring,
\[
\cO(T,\ZZ)=\ZZ p\oplus\ZZ q,
\qquad
p^2=q,
\qquad pq=qp=q^2=0.
\]
The regular action satisfies
\[
p(e)=f,
\quad p(f)=0,
\quad q(T)=0.
\]
Therefore
\[
\Def(T,\ZZ)=\UDef(T,\ZZ)=\ZZ q,
\quad
\cC(T,\ZZ)=\ZZ q.
\]
Moreover,
\begin{equation}\label{eq:nil-rings}
\cA(T,\ZZ)\cong\ZZ[p]/(p^3),
\qquad
\cR_T\cong\ZZ[p]/(p^2),
\qquad
K=(p^2).
\end{equation}
The regular extension class $\Omega(T,\ZZ)$ is nonzero.
\end{theorem}

\begin{proof}
This is the rank-one case of \cref{thm:tensor-structure}. Explicitly, in $T\otimes T$ the middle relation identifies $f\otimes e$ and $e\otimes f$ and kills $f\otimes f$. Taking $p=e\otimes e$ and $q=f\otimes e=e\otimes f$ gives the multiplication table. The regular action is immediate.

Adjoining an identity gives the first ring in \eqref{eq:nil-rings}; quotienting by the regular kernel gives the second. If a ring section existed, the class of $p$ in $\ZZ[p]/(p^2)$ would have a square-zero lift. Every lift is $p+ap^2$, whose square is $p^2\ne0$. Apply \cref{thm:HHclass}.
\end{proof}

\begin{theorem}[Periodic resolution]\label{thm:nil-periodic}
As a left or right $\cA=\ZZ[p]/(p^3)$-module, $\cR=\ZZ[p]/(p^2)$ has a two-periodic free resolution
\begin{equation}\label{eq:nil-resolution}
\cdots\xrightarrow{\,p\,}\cA
\xrightarrow{\,p^2\,}\cA
\xrightarrow{\,p\,}\cA
\xrightarrow{\,p^2\,}\cA
\longrightarrow\cR\longrightarrow0.
\end{equation}
Let $Z=\cA/(p)\cong\ZZ$. Then, for every $n\ge1$,
\begin{equation}\label{eq:nil-allZ}
\Tor_n^{\cA}(\cR,Z)\cong\ZZ,
\qquad
\Ext^n_{\cA}(\cR,Z)\cong\ZZ.
\end{equation}
With coefficients in $\cR$, for every $j\ge0$,
\begin{align}
\Tor_{2j+1}^{\cA}(\cR,\cR)&\cong\cR/p\cR\cong\ZZ,\label{eq:nil-torodd}\\
\Tor_{2j+2}^{\cA}(\cR,\cR)&\cong p\cR\cong\ZZ,\label{eq:nil-toreven}\\
\Ext^{2j+1}_{\cA}(\cR,\cR)&\cong p\cR\cong\ZZ,\label{eq:nil-extodd}\\
\Ext^{2j+2}_{\cA}(\cR,\cR)&\cong\cR/p\cR\cong\ZZ.\label{eq:nil-exteven}
\end{align}
\end{theorem}

\begin{proof}
Multiplication by $p^2$ has kernel $(p)$, and multiplication by $p$ has kernel $(p^2)$. Hence \eqref{eq:nil-resolution} is exact. Tensoring with $Z$ makes every differential zero, as does applying $\Hom_{\cA}(-,Z)$, giving \eqref{eq:nil-allZ}.

After tensoring with $\cR$, multiplication by $p^2$ is zero and multiplication by $p$ has image and kernel $p\cR$. The homology gives \eqref{eq:nil-torodd}--\eqref{eq:nil-toreven}. Applying $\Hom(-,\cR)$ gives the alternating cochain complex with maps $0$ and multiplication by $p$, yielding \eqref{eq:nil-extodd}--\eqref{eq:nil-exteven}.
\end{proof}

\begin{remark}
The Poincar\'e series of the positive-degree groups in \eqref{eq:nil-allZ}, counted by free rank, is
\[
\sum_{n\ge1}\rank\Tor_n^{\cA}(\cR,Z)t^n=\frac{t}{1-t}.
\]
Thus the invisible operator $p^2$ generates a homological signal in every positive degree, not merely in $\Tor_1$.
\end{remark}

\section{A nonunital rank-one family with an infinite defect tower}\label{sec:rank-one}

Fix $d\ge2$. Let $T_d=d\ZZ$, let $\Gamma=\ZZ$, and set
\[
x\,n\,y\,m\,z=nmxyz.
\]
This is induced by the nonunital rng $d\ZZ$.

\begin{theorem}\label{thm:Td-envelope}
Let $x=d$ be the additive generator of $T_d$. Then
\[
\cO(T_d,\ZZ)=\ZZ t,
\qquad t=x\otimes x,
\qquad t^2=d^2t.
\]
The intrinsic regular representation is injective and sends $t$ to $d^2\id_{T_d}$. Hence
\[
\Def(T_d,\ZZ)=0.
\]
For
\[
\cA_d=\ZZ[t]/(t^2-d^2t),
\qquad
\cR_d\cong\ZZ
\]
with $t\mapsto d^2$, the unital defect is generated by
\[
\kappa=t-d^2.
\]
It satisfies
\[
\kappa^2=-d^2\kappa,
\qquad
K^n=d^{2(n-1)}K.
\]
Consequently,
\begin{equation}\label{eq:Td-layers}
\cC_n(T_d,\ZZ)=K^n/K^{n+1}\cong\ZZ/d^2\ZZ
\quad\text{for every }n\ge1.
\end{equation}
\end{theorem}

\begin{proof}
The rank-one tensor group is free on $t=x\otimes x$, the middle relation is tautological, and
\[
t^2=x^3\otimes x=d^3\otimes d=d^2t.
\]
The regular action is multiplication by $d^2$, hence injective. The remaining assertions follow by substituting $t=\kappa+d^2$ into $t^2=d^2t$.
\end{proof}

\begin{theorem}[Alternating periodicity]\label{thm:Td-periodic}
As a $\cA_d$-module, $\cR_d=\cA_d/(\kappa)$ has the two-periodic free resolution
\begin{equation}\label{eq:Td-resolution}
\cdots\xrightarrow{\,t\,}\cA_d
\xrightarrow{\,\kappa\,}\cA_d
\xrightarrow{\,t\,}\cA_d
\xrightarrow{\,\kappa\,}\cA_d
\longrightarrow\cR_d\longrightarrow0.
\end{equation}
For every $j\ge0$,
\begin{align}
\Tor_{2j+1}^{\cA_d}(\cR_d,\cR_d)&\cong\ZZ/d^2\ZZ,
&\Tor_{2j+2}^{\cA_d}(\cR_d,\cR_d)&=0,\label{eq:Td-tor}\\
\Ext^{2j+1}_{\cA_d}(\cR_d,\cR_d)&=0,
&\Ext^{2j+2}_{\cA_d}(\cR_d,\cR_d)&\cong\ZZ/d^2\ZZ.\label{eq:Td-ext}
\end{align}
\end{theorem}

\begin{proof}
One has $\ann(\kappa)=(t)$ and $\ann(t)=(\kappa)$, proving exactness of \eqref{eq:Td-resolution}. On $\cR_d\cong\ZZ$, the element $\kappa$ acts as zero and $t$ acts as multiplication by $d^2$. The homology and cohomology of the resulting alternating complexes give \eqref{eq:Td-tor} and \eqref{eq:Td-ext}.
\end{proof}

\begin{remark}[Matrix-factorization position]
The maps in \eqref{eq:Td-resolution} form the rank-one matrix factorization
$(t,t-d^2)$ of $t(t-d^2)$ over $\ZZ[t]$ \cite{Eisenbud1980}.  Thus the
two-periodicity itself is the classical hypersurface mechanism.  What is
specific here is that the nonunital ternary operator construction produces
this factorization while its intrinsic defect is zero and its unitized defect
has the infinite tower \eqref{eq:Td-layers}.
\end{remark}

\begin{remark}
The intrinsic defect vanishes, but the unital defect has infinite depth and identical cyclic graded layers. This shows that intrinsic slot invisibility and identity normalization are genuinely different mechanisms. The degree-support generating series of the nonzero groups are
\[
\sum_{j\ge0}t^{2j+1}=\frac{t}{1-t^2},
\qquad
\sum_{j\ge0}t^{2j+2}=\frac{t^2}{1-t^2}.
\]
\end{remark}

\section{A torsion defect on a torsion-free ternary group}\label{sec:sign}

Let $E=\ZZ e_1\oplus\ZZ e_2$, $W=\ZZ f$, and define a surjective tensor $\mu:E^{\otimes3}\to W$ by
\begin{equation}\label{eq:sign-tensor}
\mu(e_i,e_j,e_k)=
\begin{cases}
-f,&(i,j,k)=(2,2,2),\\
f,&\text{otherwise}.
\end{cases}
\end{equation}
Let $T_{\pm}=E\oplus W$ be the corresponding two-step ternary ring.

Write $p_{ij}=e_i\otimes e_j$. With rows indexed by the final input $e_1,e_2$, the visible matrix is
\begin{equation}\label{eq:sign-Theta}
\Theta_\mu=
\begin{pmatrix}
1&1&1&1\\
1&1&1&-1
\end{pmatrix}.
\end{equation}
Hence
\[
\Ker\Theta_\mu
=\ZZ(p_{11}-p_{12})\oplus\ZZ(p_{11}-p_{21}).
\]

\begin{theorem}\label{thm:sign-defect}
For the tensor \eqref{eq:sign-tensor},
\[
Q_\mu\cong\ZZ/2\ZZ.
\]
If $q$ denotes its nonzero element, every product of two basis elements $p_{ij}$ equals $q$. Therefore
\begin{equation}\label{eq:sign-D}
\Def(T_{\pm},\ZZ)
=\UDef(T_{\pm},\ZZ)
=\cC(T_{\pm},\ZZ)
\cong\ZZ^2\oplus\ZZ/2\ZZ.
\end{equation}
Moreover, the operator rng annihilates the defect on both sides, and the regular square-zero extension is nonsplit.
\end{theorem}

\begin{proof}
Let $L_i=f\otimes e_i$ and $R_i=e_i\otimes f$. The relations are
\[
\mu(e_a,e_b,e_c)\otimes e_d
=e_a\otimes\mu(e_b,e_c,e_d).
\]
The choices $(a,b,c,d)=(1,1,1,1)$, $(1,1,1,2)$ and $(2,1,1,1)$ imply
\[
L_1=L_2=R_1=R_2=q.
\]
The choice $(2,2,2,1)$ gives $-L_1=R_2$, hence $2q=0$. Conversely, sending all four generators to $1\in\ZZ/2\ZZ$ respects every relation, so $q\ne0$ and $Q_\mu\cong\ZZ/2\ZZ$.

The product $p_{ab}p_{cd}$ is $\pm L_d$, which is always $q$ because $q=-q$. The kernel generators in \eqref{eq:sign-Theta} are differences of basis tensors, so multiplication by any $p_{ab}$ kills them; the defect also annihilates the operator rng on the left by \cref{thm:tensor-structure}. This proves \eqref{eq:sign-D} and the two-sided annihilation statement.

Let $v$ be the visible operator represented by $p_{11}$. Then $v^2=0$ in the regular image, while every lift $p_{11}+k$ has square $q\ne0$. Apply \cref{crit:nonsplit}.
\end{proof}

\begin{corollary}\label{cor:sign-homology}
Let $Z=\ZZ$ be the $\cR_{T_\pm}$-module on which the visible square-zero ideal acts trivially. Then
\[
\Tor_1^{\cA(T_\pm,\ZZ)}(\cR_{T_\pm},Z)
\cong\ZZ^2\oplus\ZZ/2\ZZ,
\]
while
\[
\Ext^1_{\cA(T_\pm,\ZZ)}(\cR_{T_\pm},Z)
\cong\ZZ^2.
\]
\end{corollary}

\begin{proof}
The regular image acts on both the defect and $Z$ through its integer quotient. Apply \cref{thm:conormal}; $\Hom_{\ZZ}(\ZZ/2,\ZZ)=0$.
\end{proof}

\begin{proposition}[Characteristic-$2$ jump]\label{prop:sign-jump}
For every odd prime $p$,
\[
\dim_{\FF_p}\Def_{\FF_p}(T_{\pm}\otimes\FF_p)=2.
\]
In characteristic $2$,
\[
\dim_{\FF_2}\Def_{\FF_2}(T_{\pm}\otimes\FF_2)=4.
\]
Thus the defect jump at $2$ equals $2$.
\end{proposition}

\begin{proof}
The Smith normal form of $B_\mu$ is
\[
\operatorname{diag}(1,1,1,2).
\]
For odd $p$, \eqref{eq:sign-Theta} has rank $2$ and $B_\mu$ has rank $4$, giving dimension $(4-2)+(4-4)=2$. Modulo $2$, the two rows of \eqref{eq:sign-Theta} coincide, so its rank is $1$, while $B_\mu$ has rank $3$. The dimension is therefore $(4-1)+(4-3)=4$.
\end{proof}

\begin{remark}
The integral defect tensored with $\FF_2$ has dimension $3$, whereas the defect of the reduced ternary system has dimension $4$. One new dimension comes from the integral $2$-torsion class $q$, and the other from the rank drop of the visible-operator matrix. This gives a concrete failure of naive nonflat base change.
\end{remark}

\section{Exact census of \texorpdfstring{$2\times2\times2$}{2x2x2} sign tensors}\label{sec:census}

Let $E=\ZZ^2$, $W=\ZZ$, and restrict to tensors whose eight structure constants lie in $\{\pm1\}$.  Every such tensor is surjective.  The finite family is small enough for exhaustive exact classification but large enough to display three distinct integral defect profiles.

\subsection{Signed-incidence interpretation}
Write the structure constants as $s_{abc}\in\{\pm1\}$, with indices in $\{0,1\}$.  Let $G_s$ be the signed bipartite multigraph with vertices
\[
L_0,L_1,R_0,R_1.
\]
For every quadruple $(a,b,c,d)$ place an edge between $R_a$ and $L_d$ with sign
\[
\varepsilon_{abcd}=s_{abc}s_{bcd}.
\]
For every pair $(a,d)$ there are four such edges, indexed by $(b,c)$.
Consequently the underlying simple graph is $K_{2,2}$ and is connected, so
\cref{lem:signed-cokernel} applies to every tensor in the census.
After multiplying a column by the unit $s_{abc}$, the corresponding column of $B_\mu$ is
\[
L_d-\varepsilon_{abcd}R_a.
\]
Thus $B_\mu$ is a signed vertex-edge incidence matrix of $G_s$ up to unimodular column operations.

\begin{lemma}[Signed-incidence cokernel]\label{lem:signed-cokernel}
Let $G$ be a connected signed multigraph and let $H_G$ be the integral matrix whose column for an edge $uv$ of sign $\varepsilon$ is $e_u-\varepsilon e_v$.  Then
\[
\coker H_G\cong
\begin{cases}
\ZZ,&G\text{ is balanced},\\
\ZZ/2\ZZ,&G\text{ is unbalanced}.
\end{cases}
\]
Here balanced means that the product of the edge signs around every cycle is positive.
\end{lemma}

\begin{proof}
Switching signs at a vertex multiplies the corresponding row by $-1$ and hence does not change the cokernel.  Switch along a spanning tree until every tree edge is positive.  The tree columns identify all vertex generators with a single generator $z$.  If $G$ is balanced, every remaining edge is positive and contributes no further relation, giving $\ZZ z$.  If $G$ is unbalanced, some remaining edge is negative and contributes $2z=0$.  Conversely every unbalanced cycle produces such a negative edge after switching.  The tree columns together with one negative edge contain a full-rank minor of determinant $\pm2$, so no additional torsion occurs.  This is the incidence-matrix form of signed-graph balance and switching \cite{Zaslavsky1982}.
\end{proof}

\begin{corollary}\label{cor:sign-Q}
For every $2\times2\times2$ sign tensor,
\[
Q_\mu\cong\ZZ\quad\text{or}\quad Q_\mu\cong\ZZ/2\ZZ.
\]
The first case occurs exactly when the signed multigraph $G_s$ is balanced.
\end{corollary}

\begin{theorem}[Complete sign-tensor census]\label{thm:sign-census}
Among the $256$ sign tensors, the triples
\[
(\rank_{\QQ}\Theta_\mu,
 \rank_{\ZZ}Q_\mu,
 \tors Q_\mu)
\]
occur with the following multiplicities:
here $\rank_{\ZZ}Q_\mu$ means the free rank of the finitely generated abelian
group $Q_\mu$.
\begin{center}
\begin{tabular}{@{}cccc@{}}
\toprule
$\rank\Theta_\mu$ & $Q_\mu$ & $\Def(T_\mu,\ZZ)$ & number \\ \midrule
$2$ & $\ZZ/2\ZZ$ & $\ZZ^2\oplus\ZZ/2\ZZ$ & $224$ \\
$1$ & $\ZZ/2\ZZ$ & $\ZZ^3\oplus\ZZ/2\ZZ$ & $28$ \\
$1$ & $\ZZ$ & $\ZZ^4$ & $4$ \\
\bottomrule
\end{tabular}
\end{center}
The four torsion-free exceptional tensors are, up to a global sign,
\[
\mu_{abc}=1
\qquad\text{and}\qquad
\mu_{abc}=(-1)^{a+b+c}
\quad(a,b,c\in\{0,1\}).
\]
Every remaining sign tensor has a nonzero $2$-primary middle-relation defect.
\end{theorem}

\begin{proof}
A $2\times4$ sign flattening has rank one precisely when its four columns lie on one of the two rational sign lines spanned by $(1,1)^t$ or $(1,-1)^t$.  Choosing the line and independently choosing the sign of each column gives $2\cdot2^4=32$ rank-one tensors; the other $224$ have rank two.

By \cref{cor:sign-Q}, it remains to count balanced signed multigraphs.  Put $s_{abc}=(-1)^{x_{abc}}$ with $x_{abc}\in\FF_2$.  Balance says that the exponents
\[
x_{abc}+x_{bcd}
\]
are independent of $(b,c)$ for each pair $(a,d)$ and have zero sum around the unique four-cycle of the underlying $K_{2,2}$.  Row reduction of these linear equations gives
\[
x_{010}=x_{001},\quad x_{011}=x_{000},\quad
x_{100}=x_{001},\quad x_{101}=x_{000},\quad
x_{110}=x_{000},\quad x_{111}=x_{001}.
\]
Thus, with $\eta=x_{000}$ and $\lambda=x_{000}+x_{001}$,
\[
x_{abc}=\eta+\lambda(a+b+c),
\qquad \eta,\lambda\in\FF_2.
\]
Hence exactly four tensors are balanced: the constant and parity tensors and their global negatives.  They have $Q_\mu\cong\ZZ$.  The remaining $252$ tensors have $Q_\mu\cong\ZZ/2\ZZ$.  Intersecting this count with the $32$ rank-one tensors gives four torsion-free and twenty-eight torsion cases, while all $224$ rank-two tensors have $2$-torsion.  Combining with
\(
\Def(T_\mu,\ZZ)\cong\Ker\Theta_\mu\oplus Q_\mu
\)
proves the table.  For completeness, the constant tensors have proportional
flattening rows, while each parity tensor has columns
$(-1)^{a+b}(1,-1)^t$; global negation does not change rank.  Hence all four
balanced tensors are among the $32$ rank-one tensors, justifying the stated
intersection count.  The supplied exact script independently enumerates all
$256$ tensors and records the Smith factors, providing a machine-checkable
certificate of the count.
\end{proof}

\begin{corollary}
A uniformly random sign tensor has torsion in its operator defect with probability $252/256=63/64$, while visible rank is maximal with probability $224/256=7/8$.
\end{corollary}

\begin{remark}
The census shows that integral middle-relation torsion is not an isolated example: it is prevalent inside this discrete coefficient ensemble.  By contrast, torsion-free hidden operators occur only at the four maximally factorized sign patterns.
\end{remark}

\section{Binary-derived collapse and its limits}

The next statement explains why unital associative examples alone cannot reveal the preceding phenomena.

\begin{proposition}\label{prop:unital-collapse}
Let $S$ be a unital associative ring, take $\Gamma=\ZZ$, and set
\[
x\,n\,y\,m\,z=nmxyz.
\]
Then $\cO(S,\ZZ)\cong S$ as rngs and $\Def(S,\ZZ)=0$.
\end{proposition}

\begin{proof}
The raw tensor group is $S\otimes S$. Multiplication $a\otimes b\mapsto ab$ descends to $\cO(S,\ZZ)$. The inverse sends $s$ to the class of $s\otimes1$. Indeed, the middle relation with the remaining entries equal to $1$ gives $a\otimes b=ab\otimes1$. The regular representation becomes the faithful left regular representation of $S$.
\end{proof}

\begin{proposition}\label{prop:unital-first-vanish}
In the setting of \cref{prop:unital-collapse}, the conormal defect acts trivially after passage to any unital $S$-module, and the first relative groups in \cref{thm:conormal} vanish.
\end{proposition}

\begin{proof}
The Dorroh kernel is generated by $(1,-1_S)$. The element $(0,1_S)$ acts as the identity on every unital $S$-module but annihilates the kernel on both sides. The tensor and Hom expressions in \cref{thm:conormal} therefore vanish.
\end{proof}

\begin{remark}
The rank-one family of \cref{sec:rank-one} shows that removing the binary identity changes the picture: the intrinsic representation remains faithful, but adjoining an identity creates an infinite arithmetic defect tower. The non-binary example of \cref{sec:nil-example} goes further by producing an intrinsic square-zero kernel and nonzero homology in every degree.
\end{remark}

\section{Exact computation and reproducibility}\label{sec:computation}

The supplementary notebook \texttt{operator\_defect\_colab.ipynb} implements the following exact workflow.

\begin{enumerate}[label=\textup{(\arabic*)},leftmargin=2.2em]
\item Input the integer tensor coefficients $\mu^u_{abc}$.
\item Build $\Theta_\mu$ and $B_\mu$ using \eqref{eq:Theta-matrix} and \eqref{eq:B-matrix}.
\item Compute rational ranks, modular ranks, Smith invariant factors and unimodular Smith change-of-basis matrices.
\item Recover the free and torsion profiles of $\Def(T_\mu,\ZZ)$.
\item Use the left Smith transformation to reduce each universal product vector to canonical torsion and free coordinates in $\coker B_\mu$, thereby obtaining the curvature table.
\item Verify the complementary Smith profiles of the visible and middle atomic families for selected parameters.
\item Verify the nilpotent and sign-tensor examples, both rank-one periodic complexes, the base-change correction term and the characteristic-jump calculation.
\item Exhaust all $256$ sign tensors, group them by exact Smith profile, and reproduce \cref{thm:sign-census}.
\item Export separate machine-readable certificates for the atomic families and the exhaustive sign census, containing ranks, invariant factors and modular dimensions for independent checking.
\end{enumerate}

The supplied executed notebook runs every listed routine with SymPy 1.14.0,
which is pinned in \texttt{requirements.txt}.  All computations use exact
integers and finite fields. There is no floating-point approximation and no
machine-learning model. The notebook is therefore a reproducibility supplement
and an exploratory tool, not a substitute for proof.

\begin{proposition}\label{prop:algorithm-correct}
For finite-free two-step tensor systems, the notebook output for the operator-defect group and its bad-prime set is mathematically complete: it is determined by the Smith normal form and the ranks of the two presentation matrices.
\end{proposition}

\begin{proof}
This is exactly \cref{thm:SNF,prop:char-jump}.
\end{proof}

\section{Consequences and open structural problems}

The preceding results separate several mathematically distinct layers of the
universal-to-regular comparison.

\subsection{Independent arithmetic design}
The realization theorem shows that the two sources of bad reduction are not coupled: visible nonsaturation and hidden middle-relation torsion may be assigned unrelated invariant factors, even though both matrices are constructed from the same trilinear tensor.  A sharper problem is now quantitative rather than existential: determine the smallest pair $(\rank E,\rank W)$ realizing a prescribed ordered pair $(G,H)$, and classify the minimal tensors up to integral change of bases.

\subsection{Finite-depth detection versus idempotent normalization}
By \cref{thm:conormal-rigidity}, first relative homology detects every nilpotent regular defect.  The only way a nonzero kernel can escape $K/K^2$ is through non-nilpotent behavior such as $K=K^2$.  This dichotomy suggests separating intrinsic finite-depth defects from idempotent normalization kernels before any homological invariant is interpreted.

\subsection{Defect syzygies}
The modules $\cS_n$ are relations among invisible operators, not merely higher powers of an ideal.  A natural next invariant is a minimal graded resolution of $\gr_K(\cA)$ over the tensor algebra $T_{\cR}(\cC)$.  Its graded Betti numbers would measure the complexity of operator identities beyond the first conormal layer and distinguish square-zero collapse from tensor-free towers.

\subsection{Arithmetic and geometric moduli}
The determinantal strata of \cref{thm:defect-strata} form a coarse moduli decomposition before quotienting by $\mathrm{GL}(E)\times\mathrm{GL}(W)$.  Their intersections with orbit closures may detect degenerations in which visible operators disappear while hidden middle relations acquire torsion.  Over $\ZZ$, the same equations define arithmetic families whose special fibers are governed by $\Delta_{\mathrm{op}}$ and the exact sequence \eqref{eq:basechange-exact}.

\subsection{Comparison with intrinsic ternary cohomology}
Carlsson and operadic cohomology deform the ternary multiplication itself, whereas $\Omega(T,\Gamma)$ deforms the extension from visible to universal operators.  Constructing a natural comparison morphism between these theories, when one exists, would identify which infinitesimal deformations of the ternary product preserve regular operator visibility and which create new hidden operators.

\subsection{Higher relative and cyclic invariants}
The periodic examples show that the quotient $\cA\to\cR_T$ carries information in every homological degree.  Change-of-rings spectral sequences, relative Hochschild homology and cyclic homology are plausible next layers, but they should be developed only after identifying ternary-specific edge maps rather than importing formal machinery without new consequences.

\begin{question}
Which graded $\cR$-algebras occur as $\gr_K\cA(T,\Gamma)$, and which occur for finite-free two-step tensor systems?
\end{question}

\begin{question}
For fixed $(r,s)$, what are the irreducible components and generic values of the determinantal defect strata in $\Hom(E^{\otimes3},W)$?
\end{question}

\begin{question}
What are the minimal ranks of $E$ and $W$ needed to realize a prescribed pair $(G,H)$, and when are the bounds in \cref{thm:independent-realization} optimal?
\end{question}

\begin{question}
When does the Hochschild curvature $\Omega(T,\Gamma)$ arise as the image of a class in intrinsic associative-triple cohomology?
\end{question}

\section{Conclusion}

The universal operator envelope of an additive ternary $\Gamma$-ring is only the starting point.  For its regular quotient we organize the standard intrinsic and unital kernels, conormal module, higher relation modules, filtered defect tower and square-zero Hochschild class.  The ternary-specific input is the canonical regular kernel and its explicit behavior in the families studied here.  Conormal rigidity shows that the first layer detects every finite-depth defect, while the relation modules distinguish square-zero collapse from tensor-free infinite towers.

For two-step tensor systems, two explicit matrices give the structure of the envelope and separate visible nonsaturation from hidden middle relations.  The atomic separation and orthogonal-sum theorems show that these arithmetic mechanisms are genuinely independent: every ordered pair of finite abelian groups occurs as the two torsion cokernels of a torsion-free two-step ternary ring.  To the authors' knowledge, this coupled realization result is new; the fact that both matrices must arise from one trilinear tensor is what prevents its reduction to arbitrary matrix-cokernel realization.

The universal envelope commutes with arbitrary scalar extension, so the nonflat phenomenon is localized in the regular kernel rather than in the presenting algebra.  The two-matrix decomposition identifies the standard Tor correction belonging to visible operators and separates it from inherited middle-relation torsion.  Maximal-rank minors and Fitting ideals give the corresponding standard determinantal strata.  Applying signed-graph switching theory to $B_\mu$ produces the exact sign-tensor census; middle-relation torsion is prevalent in this specified $256$-element ensemble.

Although the homology is computed by ordinary $\Ext$ and $\Tor$, its
ternary-specific content comes from the canonical gap between a universal
slot-sensitive operator algebra and the operators visible on the regular
ternary module, together with the filtration, deformation, arithmetic and
geometric structures carried by that gap.

\end{document}